\documentclass{article}

\usepackage[english]{babel}
\usepackage{amsfonts}
\usepackage[utf8x]{inputenc}
\usepackage[T1]{fontenc}
\usepackage{cite}
\usepackage{fullpage}
\usepackage{graphicx}
\usepackage{caption}
\usepackage{subcaption}
\usepackage{booktabs} 
\usepackage{multirow}

\usepackage{amsmath,amssymb,amsthm}
\usepackage{mathrsfs}
\numberwithin{equation}{section}
\usepackage[colorlinks=true, allcolors=blue]{hyperref}

\newtheorem{theorem}{Theorem}[section]

\newtheorem{proposition}[theorem]{Proposition}

\newtheorem{remark}[theorem]{Remark}

\newcommand{\I}{\mathcal{I}}
\newcommand{\J}{\mathcal{J}}
\newcommand{\Ir}{\widetilde{\I}}
\newcommand{\Is}{\hat{\I}}

\renewcommand{\S}{\mathcal{S}}
\newcommand{\B}{\mathcal{B}}
\newcommand{\R}{\mathcal{R}}

\newcommand{\A}{{A}}
\renewcommand{\AE}{\widetilde{\A}}
\newcommand{\E}{{M}}
\newcommand{\F}{{F}}
\newcommand{\G}{{G}}
\newcommand{\K}{{K}}
\renewcommand{\L}{{L}}
\newcommand{\T}{{T}}
\newcommand{\Z}{{Z}}
\newcommand{\C}{{C}}

\newcommand{\tB}{\widetilde{B}}

\newcommand{\defn}[1]{\emph{#1}}

\title{Recursively Preconditioned Hierarchical Interpolative Factorization for Elliptic Partial Differential Equations}
\author{Jordi Feliu-Fab\`{a}\thanks{Institute for
    Computational and Mathematical Engineering, Stanford University,
    Stanford, CA 94305, email: {\tt jfeliu@stanford.edu}}, Kenneth L. Ho\thanks{San Francisco, CA
, email: {\tt klho@alumni.caltech.edu}}, Lexing Ying\thanks{Department of Mathematics and Institute for
    Computational and Mathematical Engineering, Stanford University,
    Stanford, CA 94305, email: {\tt lexing@stanford.edu}}}
\date{}

\begin{document}

\maketitle

\begin{abstract}\label{abstract}
The hierarchical interpolative factorization for elliptic partial differential equations is a fast
algorithm for approximate sparse matrix inversion in linear or quasilinear time. Its accuracy can
degrade, however, when applied to strongly ill-conditioned problems. Here, we propose a simple
modification that can significantly improve the accuracy at no additional asymptotic cost: applying
a block Jacobi preconditioner before each level of skeletonization. This dramatically limits the
impact of the underlying system conditioning and enables the construction of robust and highly
efficient preconditioners even at quite modest compression tolerances.  Numerical examples demonstrate the performance of the new approach.
\end{abstract}

\noindent {\bf Keywords:} Recursive preconditioning; hierarchical interpolative factorization.

\section{Introduction}\label{introduction}
In this paper, we consider elliptic partial differential equations (PDE) of the form
\begin{equation}
  \label{eq1}
  -\nabla\cdot(a(x)\nabla u(x))+b(x)u(x)=f(x),\quad
  x\in\Omega\subset\mathbb{R}^{d}
\end{equation}
in two (2D, $d = 2$) and three dimensions (3D, $d = 3$), with appropriate boundary conditions on
$\partial\Omega$. Here, $a(x)$, $b(x)$, and $f(x)$ are given functions and $u(x)$ is the unknown
field. Such equations are of great importance in science and engineering and can model a wide
variety of physical phenomena. In a typical numerical solution, \eqref{eq1} is often discretized
with local schemes such as finite differences or finite elements, thus leading to a linear system
\begin{equation}
  \label{eq2}
  \A u = f,
\end{equation}
where $\A\in\mathbb{R}^{N\times N}$ is sparse and $N$ is the number of degrees of freedom (DOFs) in the discretization. Furthermore, it is common in practice for $A$ to be symmetric positive definite (SPD); we hereafter assume this to be the case. However, for non-SPD matrices, our approach can be used replacing Cholesky factorizations by $LDL^T$ or $LU$ factorizations in Section \ref{preliminaries} and Section \ref{Algorithm}.

\subsection{Background}
A significant part of research in scientific computing has been devoted to the numerical solution of \eqref{eq2}. Previous methods for its solution can largely be classified into several groups as follows.

The first group consists of classical direct methods such as Gaussian elimination or other standard matrix factorizations \cite{Golub}, nominally with $O(N^{3})$ complexity. This can be accelerated by exploiting sparsity, for instance using nested dissection (ND) \cite{Dissection} to $O(N^{3/2})$ in 2D and $O(N^2)$ in 3D for regular grids. These can still be quite prohibitive, especially for large-scale 3D problems. Therefore, although robust and highly accurate, such methods are generally not used beyond moderate problem sizes (at least in non-parallel environments).

The second group includes iterative methods \cite{Saad} such as conjugate gradient (CG) \cite{Hestenes} and multigrid \cite{Brandt,Hackbusch:multigrid}, which can achieve $O(N)$ complexity if only a small number of iterations are required. However, when $A$ is ill-conditioned or otherwise has scattered eigenvalues, as is common for high-contrast or non-smooth coefficients, the number of iterations needed can be very large. In such cases, we can expect fast convergence only if we use a good preconditioner, which can itself be a challenge to find. Furthermore, iterative methods can be inefficient for systems involving multiple right-hand sides, which are prevalent in many applications.

Bridging the two previous groups are the more recent rank-structured direct solvers, which are based on the low-rank compression of certain submatrices encountered during the factorization process. Many of these are essentially accelerated ND schemes, powered by techniques for structured dense linear algebra, and as such are associated with various dense classifications based on admissibility conditions and the use of nested bases,  including $\mathcal{H}$-matrices \cite{Grasedyck,Schmitz:2012,Schmitz:2014}, hierarchically semiseparable (HSS) matrices \cite{Xia:2009,Xia:2013:SIMAX,Xia:2013:SISC}, and hierarchically off-diagonal low-rank (HODLR) matrices \cite{Aminfar}, among others \cite{Amestoy,Gillman,Hao,HIFDE,Martinsson:2009,Pouransari,Sushnikova}. Importantly, these algorithms can be much faster than classical direct methods, with some even attaining $O(N)$ or $O(N\log N)$ complexity. Moreover, they inherently offer a speed-accuracy trade-off through the compression tolerance that is naturally suited to constructing general-purpose preconditioners. Combined with standard iterative methods, such preconditioners can enable fast and robust convergence at a far lower total cost than either a full direct solve (i.e., at high accuracy) or an unpreconditioned iterative solve.

Although rank-structured solvers have proven quite successful, they nevertheless can suffer greatly from ill-conditioning. In particular, the associated low-rank compression has traditionally been performed with respect to the forward operator so that $\A$ is well-approximated but $\A^{-1}$ may not be. Indeed, for an SPD approximation $\F = \G \G^T$ of $\A$ with $\| \A - \F \| / \| \A \| = O(\epsilon)$, it has often been observed that $\| I - \G^{-1} \A \G^{-T} \| = O(\epsilon \kappa(\A))$, where $\kappa (\cdot)$ is the condition number. Thus, ill-conditioned problems can necessitate a much higher compression accuracy than would otherwise be required in order to achieve a given solve error. This can be a significant impediment, especially for constructing low-accuracy preconditioners.

Recently, there has been some work \cite{Agullo,Xia:2017,Xing} aimed at improving the solve error and hence enabling the use of much looser compression tolerances, though only for the structured dense case so far. We highlight in particular the work of Agullo et al.\ \cite{Agullo}, which introduced a block Jacobi rescaling combined with the singular value decomposition (SVD) for low-rank approximation in the HODLR format. Specifically, they showed that (simplified and rephrased from \cite{Agullo}):

\begin{proposition}\label{prop:hodlr}
Let $A$ be a block $2 \times 2$ SPD matrix
\begin{align*}
A =
\begin{bmatrix}
A_{11} & A_{21}^T\\
A_{21} & A_{22}
\end{bmatrix} =
\begin{bmatrix}
C_1\\
& C_2
\end{bmatrix}
\begin{bmatrix}
I & B^T\\
B & I
\end{bmatrix}
\begin{bmatrix}
C_1^T\\
& C_2^T
\end{bmatrix},
\end{align*}
where each $A_{ii}$ has Cholesky decomposition $A_{ii} = C_i C_i^T$ and $B = C_2^{-1} A_{21} C_1^{-T}$. Furthermore, let $B$ have truncated SVD approximation $\tilde{B}$ (i.e., projected onto the leading singular subspace) and define
\begin{align*}
F =
\begin{bmatrix}
C_1\\
& C_2
\end{bmatrix}
\begin{bmatrix}
I & \tB^T\\
\tB & I
\end{bmatrix}
\begin{bmatrix}
C_1^T\\
& C_2^T
\end{bmatrix} = GG^T.
\end{align*}
If $\| B - \tB \| \leq \epsilon$, then $\| I - G^{-1} A G^{-T} \| \leq \epsilon$. In particular, $\kappa (G^{-1} A G^{-T}) \leq (1 + \epsilon)/(1 - \epsilon)$.
\end{proposition}

In other words, ``preconditioning'' first with block Cholesky factors before SVD compression allows very precise control over the solve error and therefore also on the convergence of preconditioned CG, which depends on the spectrum of $G^{-1} A G^{-T}$. This essentially describes an optimal one-level scheme; the multilevel extension is immediate via a recursive binary partitioning following the HODLR framework. Some minor loss of accuracy is incurred, but the overall control is still very tight. Indeed, effective preconditioners for Schur complements associated with 2D PDEs were demonstrated at tolerances up to $\epsilon \sim 0.1$. Very similar methods are presented in \cite{Xia:2017,Xing}. Other related efforts include \cite{Bebendorf:2013,Bebendorf:2016,Yang}.

However, HODLR/HSS methods have optimal linear or quasilinear complexity only if the approximation rank $\rho$ grows at most very slowly with $N$. Thus, they are really best suited for ``one-dimensional'' (1D) dense problems (by analogy with elliptic integral equations), where typically $\rho = O(\log N)$ \cite{Ho,Martinsson:2005}. This is the case for sparse PDEs in 2D \cite{Chandrasekaran}, for which elimination along, e.g., the ND ordering creates dense fill-in loaded on 1D separator edges. But it does not adequately capture the situation in 3D, where now the dense subproblems live on 2D faces, with $\rho = O(N^{1/3})$. More advanced techniques \cite{Corona,IFMM,HIFIE,Minden} are required to reduce the cost, which are not addressed by \cite{Agullo,Xia:2017,Xing}.

\subsection{Contributions}\label{sec:contrib}
Our principal goal in this paper will be to extend the ideas of \cite{Agullo} beyond the HODLR framework to a suitably general accelerated ND method capable of efficiently handling both 2D and 3D problems. We focus in particular on the hierarchical interpolative factorization (HIF) \cite{HIFDE}, a fast algorithm for computing approximate generalized Cholesky decompositions in $O(N)$ or $O(N \log N)$ time. This is achieved through alternating levels of elimination and ``skeletonization'', which eliminates DOFs from dense matrices by exploiting low-rank structure, to sparsify and reduce the dimension of the separators. For example, in 3D, elimination first reduces the problem to 2D on separator faces and then skeletonization reduces that to 1D along edges, yielding $O(N \log N)$ cost. Additionally, skeletonizing the edges themselves can further bring this down to $O(N)$, thus completing a full dimensional reduction sweep.

However Proposition \ref{prop:hodlr} is quite specific to the block $2 \times 2$ case, where the off-diagonal blocks can be simultaneously diagonalized and any error amplification suppressed due to the orthogonality of certain subspaces. This does not apply to HIF, which instead maintains a global view of all matrix blocks at each level of the factorization. Still, we can appeal to the same intuition and precondition before each round of skeletonization. The essential effect of this change can be understood heuristically as follows. First consider
\begin{align*}
\A = BB^T + E, \quad \| E \| \leq \epsilon \| \A \|,
\end{align*}
which describes a standard ``unpreconditioned'' approximation to relative precision $\epsilon$. Then the solve error is
\begin{align*}
\| I - B^{-1} \A B^{-T} \| = \| B^{-1} E B^{-T} \| \leq \epsilon \| \A \| \| B^{-1} \| \| B^{-T} \| \sim \epsilon \kappa (\A)
\end{align*}
as previously asserted. Now let
\begin{align*}
\A = C \AE C^T = C(\tB \tB^T + \tilde{E})C^T, \quad \| \tilde{E} \| \leq \epsilon \| \AE \|,
\end{align*}
i.e., the approximation is done after symmetrically preconditioning with $C$. Then $\A \approx \F = \G\G^T$ with $\G = C \tB$, so
\begin{align*}
\| I - \G^{-1} \A \G^{-T} \| = \| \tB^{-1} \tilde{E} \tB^{-T} \| \leq \epsilon \| \tB^{-1} \| \| \tB^{-T} \| \| \AE \| \sim \epsilon \kappa (\AE)
\end{align*}
and the error is now amplified only by $\kappa (\AE) \ll \kappa (\A)$. This is a much more general yet necessarily weaker result that is compatible with HIF. In particular, no detailed assumptions are made on either the block partitioning, preconditioner type, or compression method; at the same time, the accuracy is still subject to the conditioning of $\AE$, which may nevertheless be poor. The latter can be improved in the multilevel setting, where effectively $\AE$ itself is preconditioned at the next level and so on. Altogether, this outlines a ``recursively preconditioned'' HIF (PHIF) with significantly enhanced robustness to $\kappa (A)$.

We demonstrate this approach using local Cholesky factors as in \cite{Agullo} to precondition at each scale. As expected, we find substantial improvements in the solve accuracy (often by several orders of magnitude) and therefore in its effectiveness as a preconditioner, especially for ill-conditioned problems. Furthermore, the same asymptotic rank estimates are observed to hold as before and hence the computational complexity is preserved.  With this simple modification, we thus construct the first ``preconditioned'' structured nested dissection (ND) solver with optimal or near-optimal performance for PDEs in 2D and 3D.

\section{Preliminaries}\label{preliminaries}
This section reviews some key linear algebraic primitives used in PHIF. We adopt the following notation hereafter: uppercase letters ($A$, $B$, $F$, etc.)\ denote matrices; calligraphic letters ($\I$, $\J$, etc.)\ denote sets of indices, each of which is associated with a DOF; $A_{\I\J}$ is the submatrix of $A$ restricted to $\I$ and $\J$, respectively, for the rows and columns; and $[n] = \{ 1, \dots, n \}$ for $n \in \mathbb{N}$.

\subsection{Block elimination}
\label{blockelimination}
Consider an SPD matrix $\A \in \mathbb{R}^{N\times N}$ with block partitioning
\[
\A = \begin{bmatrix}
  \A_{\I\I}&\A_{\B\I}^T&\\
  \A_{\B\I}&\A_{\B\B}&\A_{\R\B}^T\\
  &\A_{\R\B}&\A_{\R\R}\\
\end{bmatrix},
\label{eqn:matrix-3block}
\]
where $[N] = \I\cup\B\cup\R$ up to permutation. This type of structure is characteristic of the sparse linear system \eqref{eq2}, which discretizes the PDE \eqref{eq1} through a partitioning of the domain $\Omega$ as a union of multiple cells with disjoint interior. For a given cell, $\I$ then represents the DOFs inside that cell, $\B$ the DOFs on the boundary, and $\R$ the remaining DOFs outside the cell and so do not interact with those in $\I$, i.e., $\A_{\R\I} = 0$.

The goal of block elimination is to zero out the $\A_{\B\I}$ block in order to decouple $\I$ from the rest. Let $\A_{\I\I}$ have Cholesky decomposition $\A_{\I\I} = \L_\I \L_\I^T$ and define the \defn{elimination matrix}
\begin{equation}
  \E_{\I}=
  \begin{bmatrix}
    \L_{\I}^{-T} & -A_{\I\I}^{-1} A_{\B\I}^{T}\\
    &I\\
    &&I
  \end{bmatrix} =
  \begin{bmatrix}
    \L_{\I}^{-T}\\
    & I\\
    && I
  \end{bmatrix}
  \begin{bmatrix}
    I & -\L_{\I}^{-1} A_{\B\I}^{T}\\
    & I\\
    && I
  \end{bmatrix} \in \mathbb{R}^{N\times N}.  \label{V_schur}
\end{equation}
Then
\[
\E_{\I}^{T}\A \E_{\I}=\begin{bmatrix}
  I&&\\
  &\AE_{\B\B}&\A_{\R\B}^T\\
  &\A_{\R\B}&\A_{\R\R}\\
\end{bmatrix}, \quad\AE_{\B\B}=\A_{\B\B}-\A_{\B\I}\A_{\I\I}^{-1}\A_{\B\I}^T.
\]
Notice that $\A_{\R\B}$ and $\A_{\R\R}$ remain unchanged while $\I$ has been fully eliminated; we are now left with a smaller problem over the restricted indices $\B \cup \R$ only.

In the context of \eqref{eq2}, we often want to perform block elimination on each of a collection of, say, $p$ cells with interior and boundary indices $\{\I_{i}\}_{i=1}^p$ and $\{\B_{i}\}_{i=1}^p$, respectively, and $\A_{\I_i, \I_j} = 0$ for all $i \neq j$. Because the Schur complement updates to each $\A_{\B_i\B_i}$ are commutative, the cells can be eliminated independently in any order, and the resulting matrix $\E=\prod_{i=1}^{p}\E_{\I_{i}}$, where each $\E_{\I_{i}}$ is given by \eqref{V_schur} with $\I=\I_{i}$, $\B=\B_{i}$, and $\R=[N]\setminus(\I_{i}\cup \B_{i})$, is well-defined. Observe then that $\E^T\A\E$ consists of the identity along $\cup_{i=1}^p \I_i$ and reduces to a subsystem in $[N]\setminus\cup_{i=1}^p\I_i$ only, i.e., the DOFs in each $\I_i$ have been eliminated.

\subsection{Skeletonization}\label{skeletonization}
Skeletonization is the generalization of block elimination for dense matrices with low-rank off-diagonal blocks. Let
\[
\A =
\begin{bmatrix}
  \A_{\I\I}&\A_{\R\I}^T\\ \A_{\R\I}&\A_{\R\R}\\
\end{bmatrix} \in \mathbb{R}^{N \times N}
\]
be SPD with $\A_{\R\I} \in\mathbb{R}^{N_\R \times N_\I}$ having numerical rank $k$ to relative precision $\epsilon$. Then from the interpolative decomposition (ID) \cite{ID}, there exists a disjoint partitioning $\I=\Ir\cup\Is$ into \defn{redundant} and \defn{skeleton} DOFs, respectively, and an \defn{interpolation matrix} $\T_{\I}\in\mathbb{R}^{k\times (N_\I-k)}$ such that
\[
\A_{\R\Ir} = \A_{\R\Is}\T_{\I} + E_{\I}, \quad \|E_{\I}\| = O(\epsilon\|\A_{\R\I}\|),
\]
i.e., any redundant column of $\A_{\R\I}$ can be well approximated by a linear combination of the skeleton columns of $\A_{\R\I}$. The choice of $\Ir$ and $\Is$ is not unique and is typically chosen so that $\|\T_\I\|$ is not too large. Following this approximation, we can rewrite $\A$ (up to permutation) as
\[
\A = \left[{\begin{array}{cc|c}
  \A_{\Ir\Ir}&\A_{\Is\Ir}^T&\A_{\R\Ir}^T\vspace*{0.1cm}\\ \vspace*{0.1cm}
  \A_{\Is\Ir}&\A_{\Is\Is}&\A_{\R\Is}^T\\ \hline
  \A_{\R \Ir}&\A_{\R \Is}&\A_{\R \R}\\
\end{array}}\right] \approx \left[{\begin{array}{cc|c}
  \A_{\Ir\Ir}&\A_{\Is\Ir}^T&\T_\I^T\A_{\R\Is}^T \vspace*{0.1cm}\\\vspace*{0.1cm}
  \A_{\Is\Ir}&\A_{\Is\Is}&\A_{\R\Is}^T\\ \hline
  \A_{\R \Is}\T_\I&\A_{\R \Is}&\A_{\R \R}\\
\end{array}}\right].
\]
Introducing the \defn{zeroing matrix} $Z_\I$ as follows and applying it on both sides leads to
\begin{equation}
  \label{eq_T}
  \Z_{\I}^{T}\A\Z_{\I}\approx\left[{\begin{array}{cc|c}
        \AE_{\Ir\Ir}&\AE_{\Is\Ir}^T& \vspace*{0.1cm}\\\vspace*{0.1cm} \AE_{\Is\Ir}&\A_{\Is\Is}&\A_{\R\Is}^T\\ \hline
        &\A_{\R\Is}&\A_{\R\R}\\
    \end{array}}\right],
  \quad
  \Z_{\I}=\left[{\begin{array}{cc|c}
        I&&\\
        -\T_{\I}&I&\\ \hline
        &&I\\
    \end{array}}\right]\in\mathbb{R}^{N\times N},
\end{equation}
where only the terms $\AE_{\Ir\Ir}$ and $\AE_{\Is\Ir}$ are updated with
\[
\AE_{\Ir\Ir}=\A_{\Ir\Ir}-\T_{\I}^{T}\A_{\Is\Ir}-\A_{\Is\Ir}^T\T_{\I}+\T_{\I}^{T}\A_{\Is\Is}\T_{\I},
\quad
\AE_{\Is\Ir}=\A_{\Is\Ir}-\A_{\Is\Is}\T_{\I}.
\]
We can now use block elimination to eliminate $\Ir$, assuming that $\AE_{\Ir\Ir}$ remains SPD. Letting $\E_{\Ir}$ be the elimination matrix so defined and introducing
$\K_{\I} \equiv \Z_{\I} \E_{\Ir}$, we therefore have
\[
\E_{\Ir}^{T}\Z_{\I}^{T}\A\Z_{\I}\E_{\Ir}\equiv\K_{\I}^{T}\A\K_{\I}
\approx\left[{\begin{array}{cc|c}
      I&&\\
      &\AE_{\Is\Is}&\A_{\R\Is}^T \vspace*{0.1cm}\\  \hline
      &\A_{\R\Is}&\A_{\R\R}\\
  \end{array}}\right], \quad \K_{\I} = \Z_{\I} \E_{\Ir}
\]
with $\AE_{\Is\Is}=\A_{\Is\Is}-\AE_{\Is\Ir}\AE_{\Ir\Ir}^{-1}\AE_{\Is\Ir}^T$. In other words, the
redundant DOFs have been eliminated without any modification to the matrix entries involving $\R$.

\begin{remark}\label{rmk:skeleton-spd}
Recall Weyl's inequality: for any perturbation matrix $P$, $|\lambda_i (A + P) - \lambda_i (A)| \leq \| P \|$, where $\lambda_i (\cdot)$ is the $i$th ordered eigenvalue. An immediate consequence is that if $\|E_{\Ir}\| < \lambda_{\min} (\A)$ in the above, where $\lambda_{\min} (\cdot)$ is the smallest eigenvalue, then $\lambda_{\min} (\Z_{\I}^T \A \Z_{\I}) > 0$, so $\AE_{\Ir\Ir}$ is SPD.
\end{remark}

As before, we often want to perform skeletonization for each of a collection of disjoint index sets
$\{\I_i\}_{i=1}^p$. But $\K_{\I_1} \K_{\I_2} \neq \K_{\I_2} \K_{\I_1}$ in general since
skeletonizing, say, $\I_1$ first can change the set $\R$ for $\I_2$. Commutativity can be restored
by instead using the ID to first find the redundant and skeleton DOFs for all $\I_i$ before zeroing
and elimination. We will not distinguish between these two approaches as they make little difference
from a practical point of view, though we note that the former ``sequential'' method is typically
faster in a serial setting, while the latter is amenable to parallelization \cite{DHIF}. In
particular, we will simply write $\K = \prod_{i=1}^p \K_{\I_i}$ for the aggregate skeletonization
operator, where the product is understood to be taken in some appropriate order. The skeletonized
matrix $\K^T \A \K$ then consists of the identity along $\cup_{i=1}^p \Ir_i$ and reduces to a
subsystem in $[N] \setminus \cup_{i=1}^p \Ir_i = \cup_{i=1}^p \Is_i$ only.

\subsection{Block Jacobi preconditioning}
\label{sec:jacobi-precond}
Let $\A$ be as in \eqref{eqn:matrix-3block} and suppose that $\A_{\I\I}$ has Cholesky decomposition $\A_{\I\I} = L_\I L_\I^T$. Then
\begin{equation}
 \C_\I^T A \C_\I =
 \begin{bmatrix}
  I & \AE_{\B\I}^T\\
  \AE_{\B\I} & \A_{\B\B} & \A_{\R\B}^T\\
  & \A_{\R\B} & \A_{\R\R}
 \end{bmatrix}, \quad \C_\I =
 \begin{bmatrix}
  L_\I^{-T}\\
  & I\\
  && I
 \end{bmatrix} \in \mathbb{R}^{N \times N},
 \label{eqn:jacobi-precond}
\end{equation}
where $\AE_{\B\I} = \A_{\B\I} L_\I^{-T}$. Since $|\B|$ is assumed small relative to $|\R|$, only a limited number of matrix entries are modified. Moreover, the \defn{preconditioning matrix} $\C_\I$ in \eqref{eqn:jacobi-precond} is block diagonal, so any collection $\{ \I_i \}_{i = 1}^p$ of disjoint index sets can be block preconditioned independently via $\C = \prod_{i = 1}^p \C_{\I_i}$; the preconditioned matrix $\C^T \A \C$ has unit block diagonal.

\section{Algorithm}\label{Algorithm}
\begin{figure}
\centering \captionsetup[subfigure]{labelformat=empty}
\begin{subfigure}{0.225\textwidth}
 \centering
 \includegraphics[width=\textwidth,trim=1.65cm 1.65cm 0.75cm 1.15cm,clip]{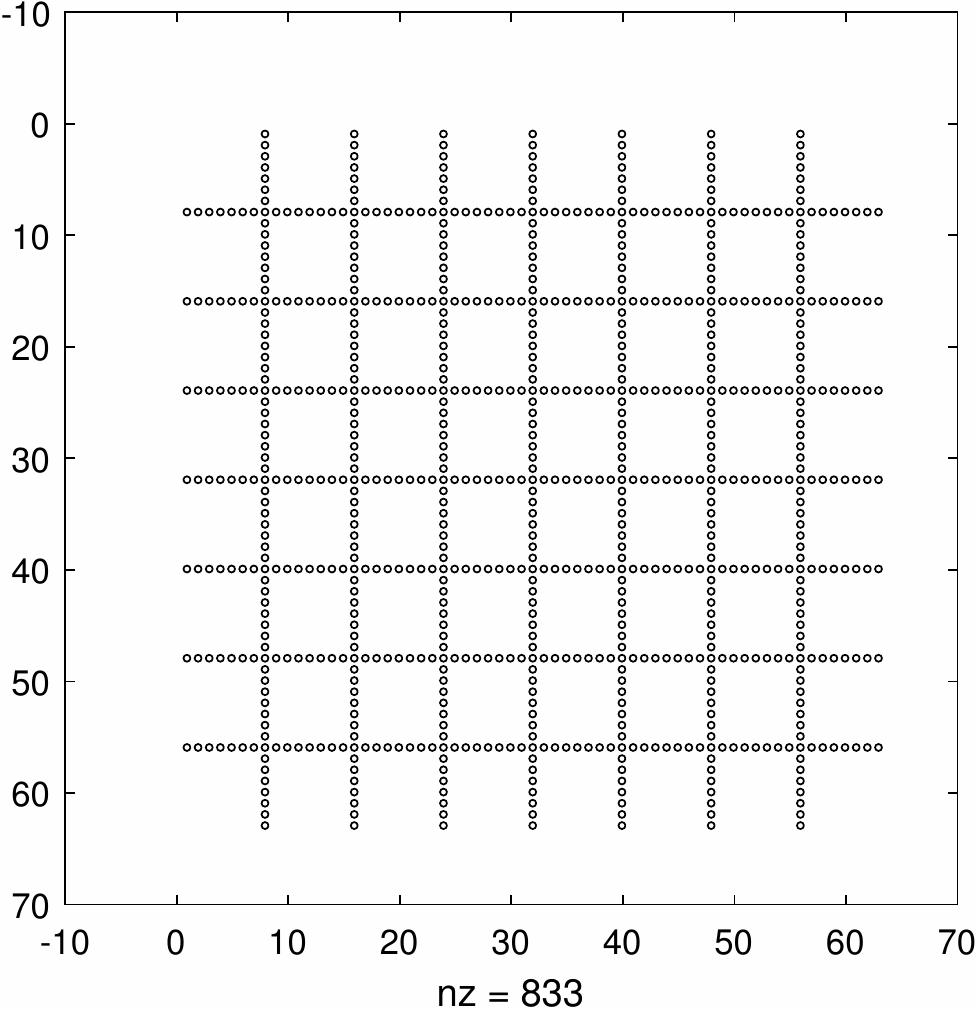}
 \caption{$\ell=0$ cell elimination}
\end{subfigure}
\begin{subfigure}{0.225\textwidth}
 \centering
 \includegraphics[width=\textwidth,trim=1.65cm 1.65cm 0.75cm 1.15cm,clip]{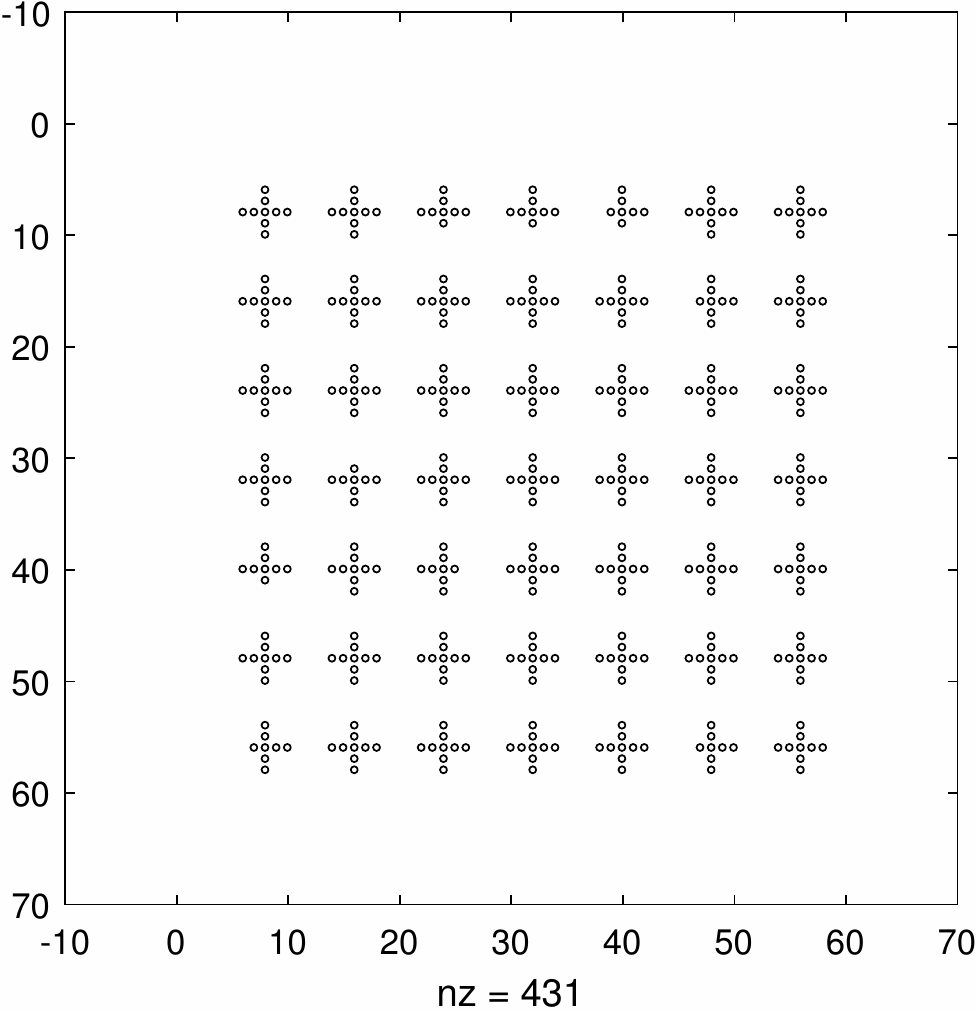}
 \caption{$\ell=0$ skeletonization}
\end{subfigure}
\begin{subfigure}{0.225\textwidth}
 \centering
 \includegraphics[width=\textwidth,trim=1.65cm 1.65cm 0.75cm 1.15cm,clip]{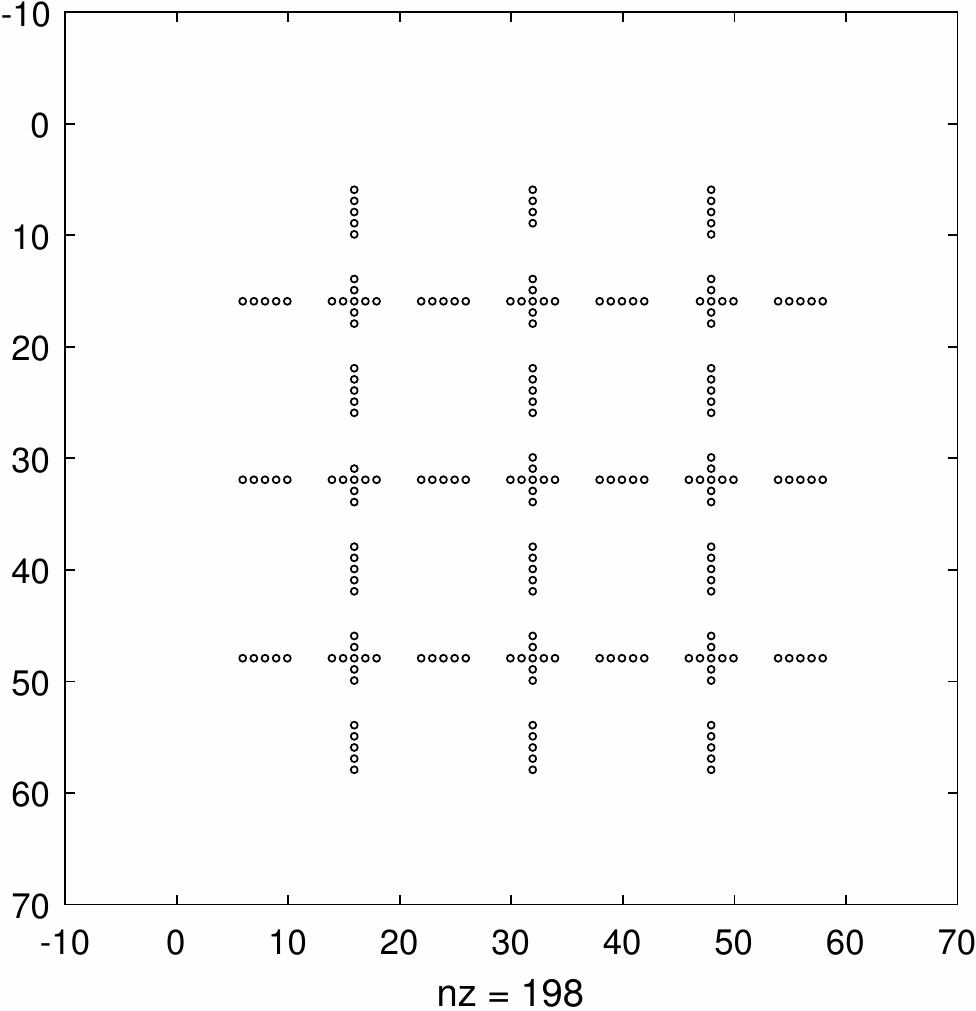}
 \caption{$\ell=1$ cell elimination}
\end{subfigure}\\
\begin{subfigure}{0.225\textwidth}
 \centering
 \includegraphics[width=\textwidth,trim=1.65cm 1.65cm 0.75cm 1.15cm,clip]{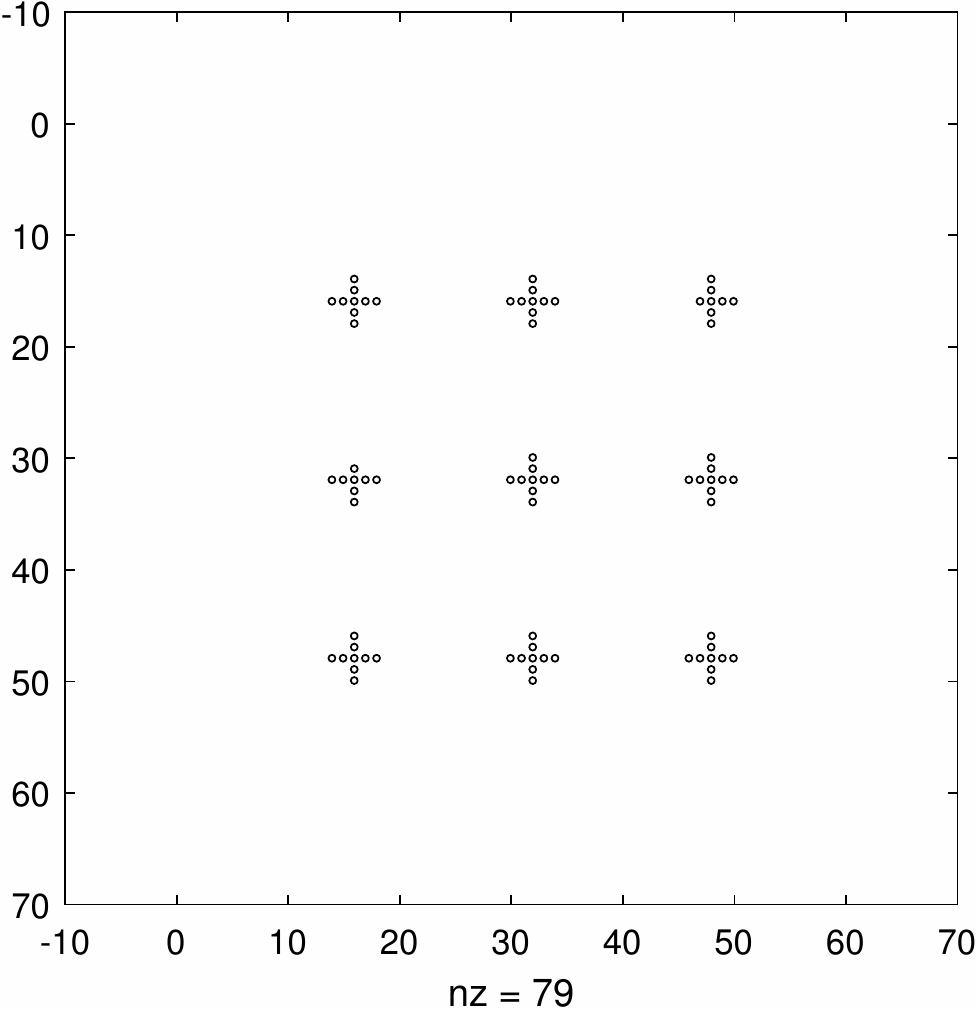}
 \caption{$\ell=1$ skeletonization}
\end{subfigure}
\begin{subfigure}{0.225\textwidth}
 \centering
 \includegraphics[width=\textwidth,trim=1.65cm 1.65cm 0.75cm 1.15cm,clip]{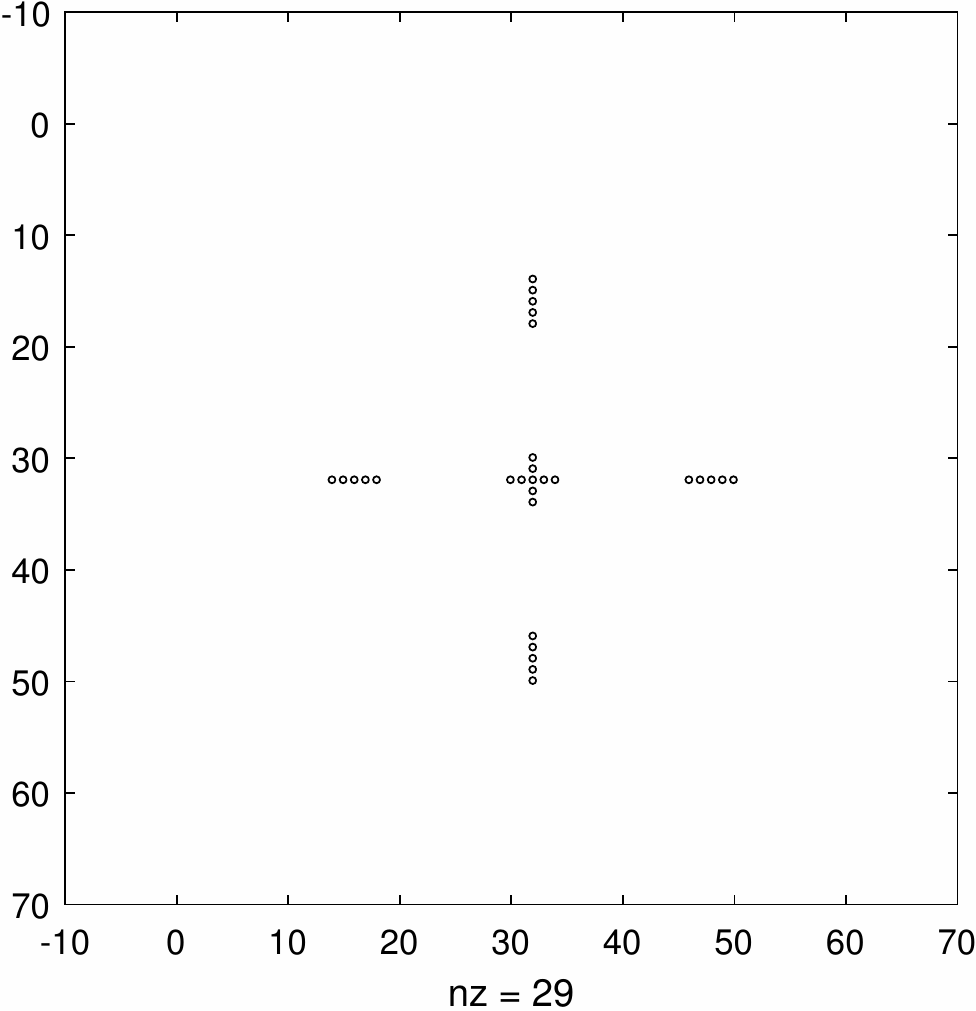}
 \caption{$\ell=2$ cell elimination}
\end{subfigure}
\begin{subfigure}{0.225\textwidth}
 \centering
 \includegraphics[width=\textwidth,trim=1.65cm 1.65cm 0.75cm 1.15cm,clip]{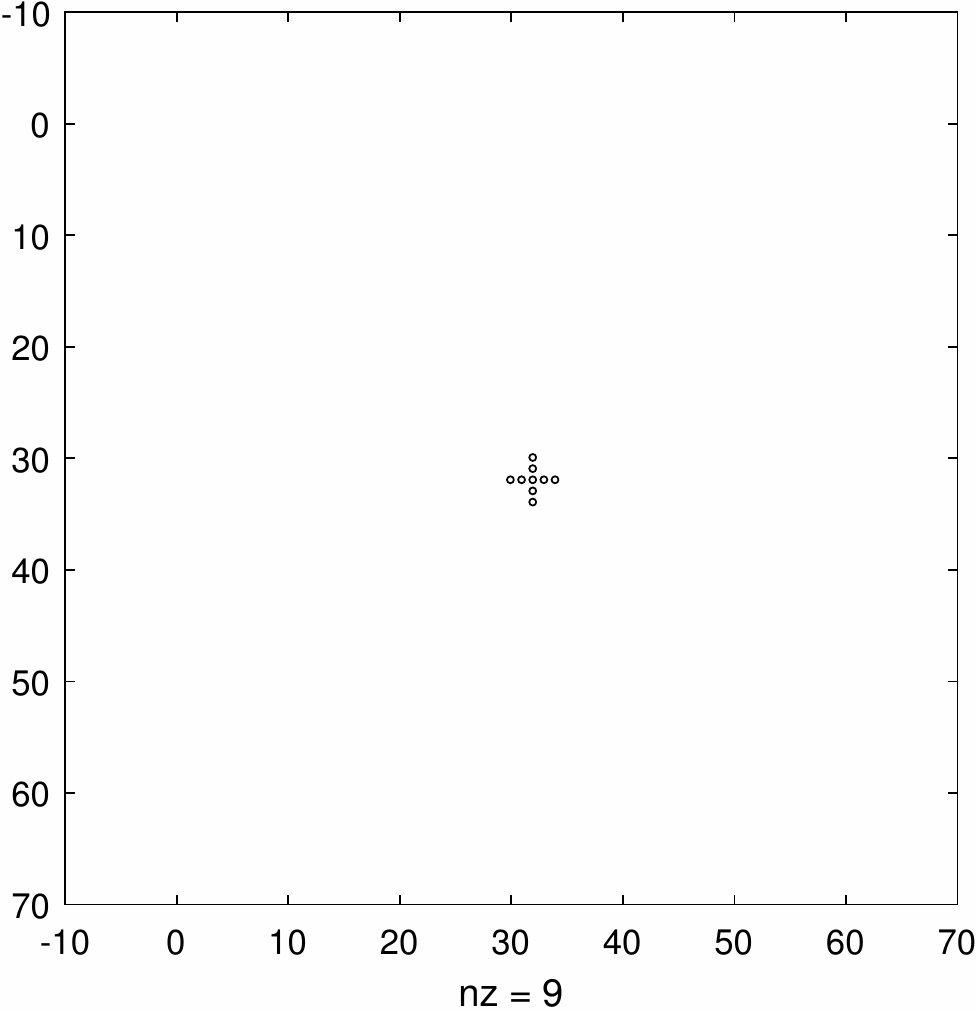}
 \caption{$\ell=2$ skeletonization}
\end{subfigure}
\caption{Active DOFs at each level $\ell$ of HIF in 2D.}
\label{HIFDE-fig}
\end{figure}

\begin{figure}
\centering \captionsetup[subfigure]{labelformat=empty}
\begin{subfigure}{0.225\textwidth}
 \centering
 \includegraphics[width=\textwidth,trim=2.1cm 2.1cm 2.1cm 1.9cm,clip]{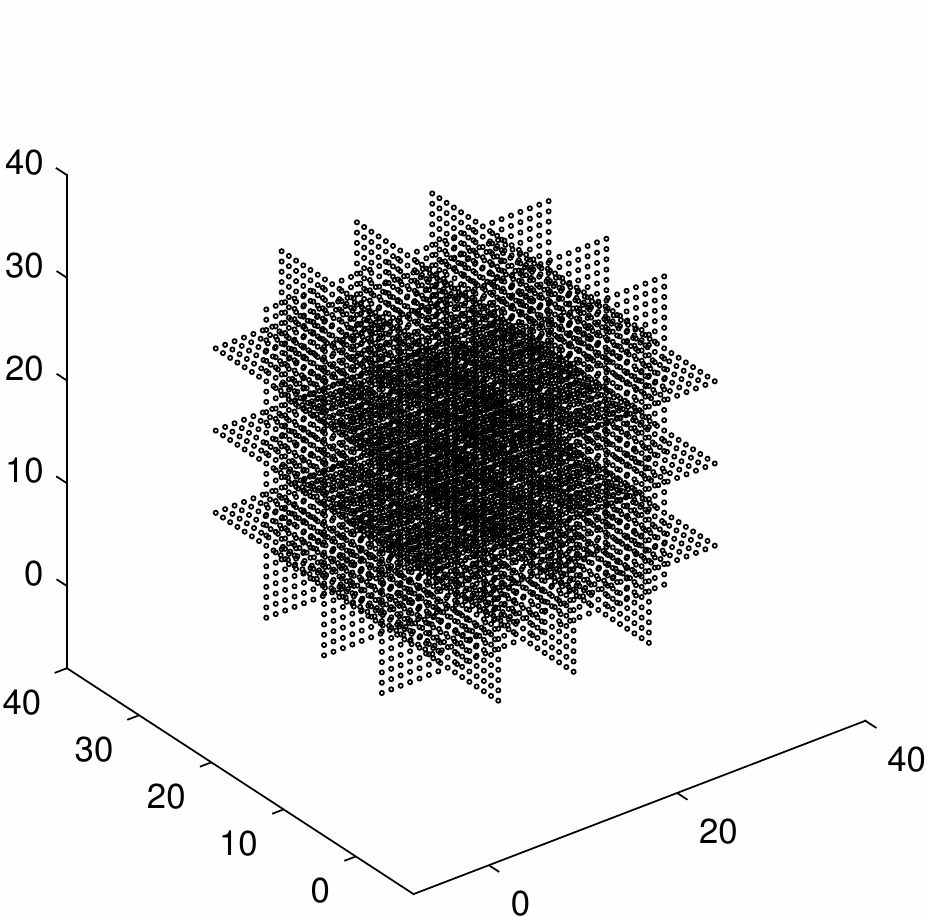}
 \caption{$\ell=0$ cell elimination}
\end{subfigure}
\begin{subfigure}{0.225\textwidth}
 \centering
 \includegraphics[width=\textwidth,trim=2.1cm 2.1cm 2.1cm 1.9cm,clip]{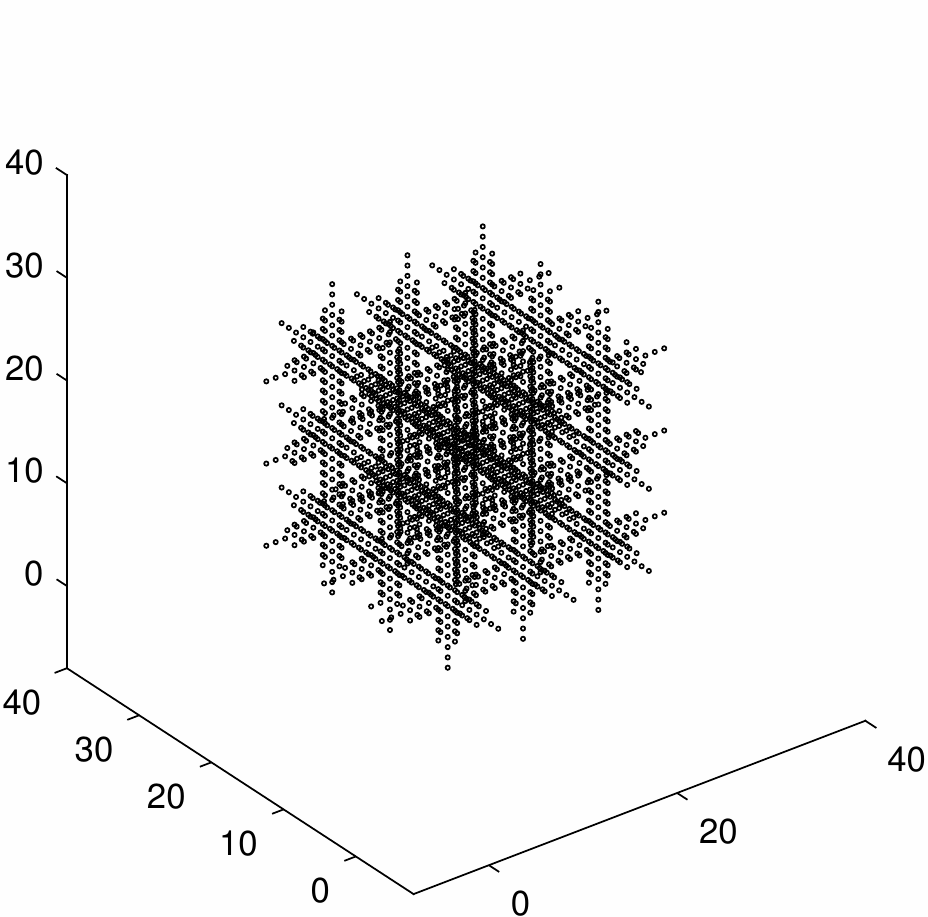}
 \caption{$\ell=0$ skeletonization}
\end{subfigure}
\begin{subfigure}{0.225\textwidth}
 \centering
 \includegraphics[width=\textwidth,trim=2.1cm 2.1cm 2.1cm 1.9cm,clip]{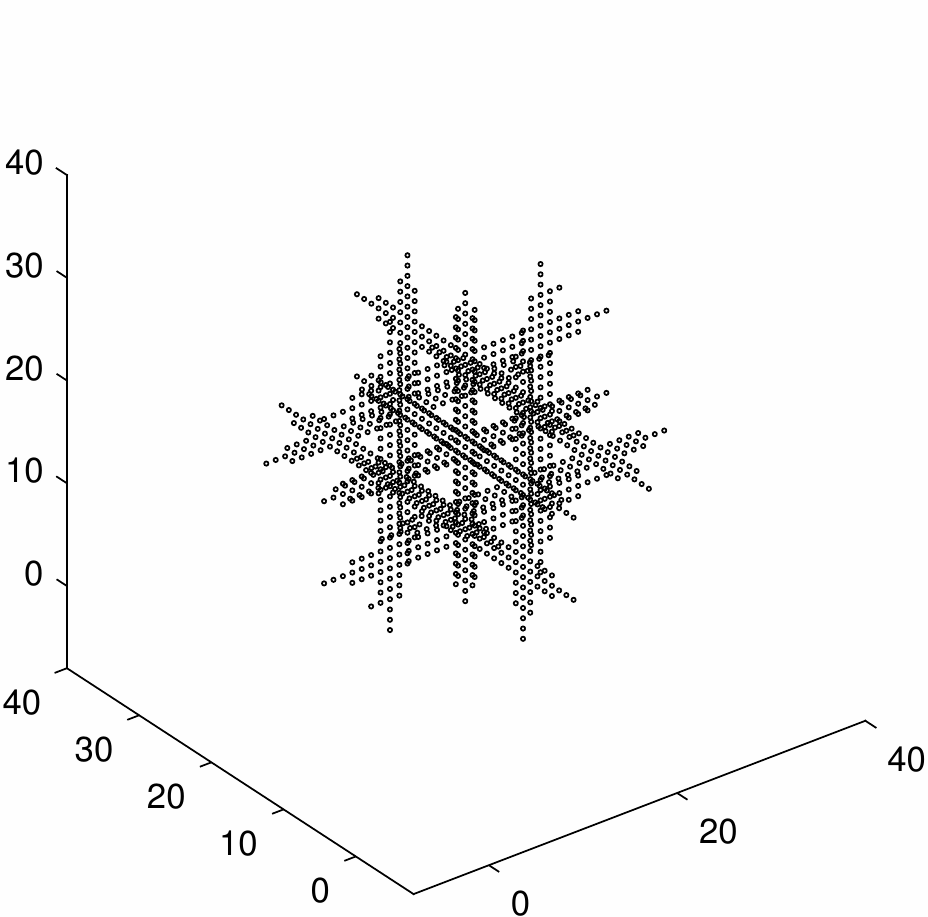}
 \caption{$\ell=1$ cell elimination}
\end{subfigure}
\begin{subfigure}{0.225\textwidth}
 \centering
 \includegraphics[width=\textwidth,trim=2.1cm 2.1cm 2.1cm 1.9cm,clip]{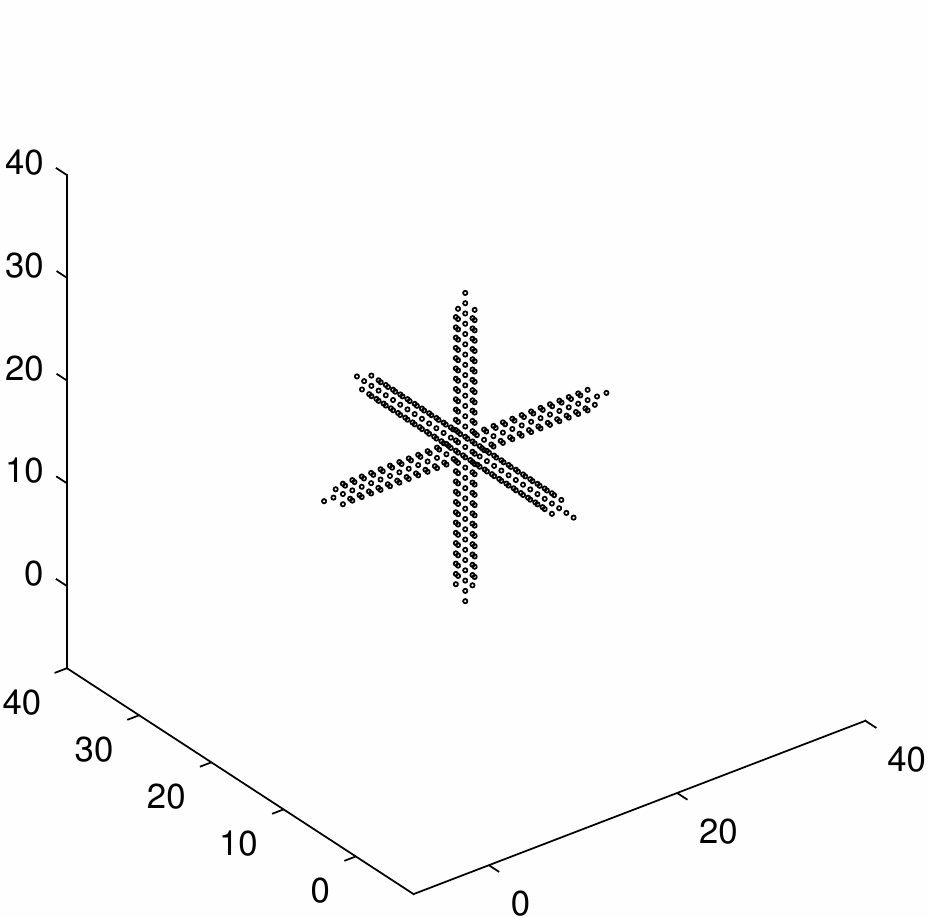}
 \caption{$\ell=1$ skeletonization }
\end{subfigure}
\caption{Active DOFs at each level $\ell$ of HIF in 3D.}
\label{HIFDE-fig3d}
\end{figure}

HIF utilizes alternating levels of block elimination and skeletonization to sparsify and eliminate from the system matrix at each problem scale, following the ND tree. It has a natural geometric interpretation, with both elimination and skeletonization acting as dimensional reduction operators. In 2D, block elimination first reduces from 2D cells to 1D edges then skeletonization reduces that to ``zero-dimensional'' points, thus achieving $O(N)$ complexity. In 3D, the same procedure gives a reduction from 3D cells to 2D faces to 1D edges, with near-optimal $O(N \log N)$ total cost; estimated $O(N)$ scaling can be recovered by additionally skeletonizing the edges, but we will not consider that here for simplicity. Figures \ref{HIFDE-fig} and \ref{HIFDE-fig3d} provide such a geometric view by showing the remaining DOFs after each step for some small examples. We refer the reader to \cite{HIFDE} for further details.

Our new PHIF follows the same framework but now adds an extra preconditioning step before each level of skeletonization. As motivated in Section \ref{sec:contrib}, this serves to control the error amplification in the matrix inverse. While the algorithm can be described in quite general terms, requiring only some sensible geometric partitioning in terms of cells, faces, and edges, as appropriate---it is robust to exactly how these are defined---we will be specific here in order to fix ideas.

\subsection{Two dimensions}
Without loss of generality, consider the PDE \eqref{eq1} on $\Omega=(0,1)^{2}$ with Dirichlet boundary conditions, discretized using finite differences via the standard five-point stencil over a uniform grid with step size $h$. We assume that the grid size $n=1/h$ in each dimension satisfies $n = 2^L m$ for integer $L$ and $m$ with $m=O(1)$. The DOFs are then the solution values $u_j = u(x_j)$ at the grid points $x_j = jh = (j_1, j_2) \cdot h$ for integers $1 \leq j_1, j_2 \leq n-1$. The resulting matrix $\A$ in \eqref{eq2} is SPD and sparse, consisting only of nearest-neighbor interactions with $\{ u_{j \pm e_i} \}_{i = 1}^2$ for each $u_j$, where $e_i$ is the $i$th unit coordinate vector. The total number of DOFs is $N = (n - 1)^{2}$.

\begin{figure}
\centering
\begin{subfigure}{0.2\textwidth}
 \centering
 \includegraphics[width=\textwidth]{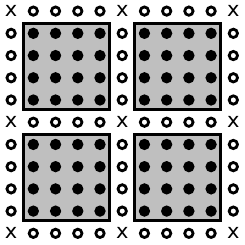}
 \caption{interior cells}
\end{subfigure}
\hspace{2em}
\begin{subfigure}{0.2\textwidth}
 \centering
 \includegraphics[width=\textwidth]{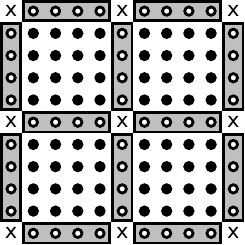}
 \caption{edges}
\end{subfigure}
\hspace{2em}
\begin{subfigure}{0.2\textwidth}
 \centering
 \includegraphics[width=\textwidth]{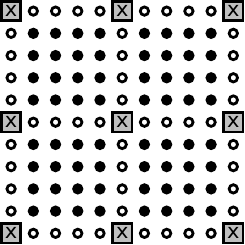}
 \caption{corners}
\end{subfigure}
\caption{Grouping of DOFs into interior cells, edges, and corners in 2D.}
\label{fig:dof-groups}
\end{figure}

Define a uniform quadtree on the domain $\Omega$ so that it partitions into $2^{(L-\ell)}\times 2^{(L-\ell)}$ square cells at each level $\ell = 0, 1, \dots, L$, going from the leaves to the root. Each such cell covers $(2^{\ell} m + 1) \times (2^{\ell} m + 1)$ grid points (or ghost points for those touching the boundary); among these,
\begin{itemize}
\item $(2^{\ell} m - 1)^2$ are \defn{interior} points;
\item $4(2^{\ell} m - 1)$ are \defn{edge} points, where each edge is shared between up to two cells; and
\item $4$ are \defn{corner} points, each shared between up to four cells.
\end{itemize}
See Figure \ref{fig:dof-groups} for a schematic. Let $p_\ell$ be the total number of cells on level $\ell$, $q_\ell \sim 2 p_\ell$ the total number of edges, and $r_\ell \sim 3 p_\ell$ the total number of edges and corners altogether.

The algorithm proceeds by eliminating DOFs level by level. Let $\S_\ell$ be the remaining \defn{active} DOFs at level $\ell$ and $\A_\ell$ the corresponding state of the matrix; initially, $\S_0 = [N]$ and $\A_0 = \A$. Then we perform the following in sequence at each level $\ell = 0, 1, \dots, L - 1$:
\begin{enumerate}
\item
 \textbf{Cell elimination.} Group the active DOFs $\S_\ell$ by interior cells and let $\{ \I_{\ell,i} \}_{i = 1}^{p_\ell}$ be the collection of all such index sets. Block elimination with respect to $\{ \I_{\ell,i} \}_{i = 1}^{p_\ell}$ as in Section \ref{blockelimination} then gives
 \[
  \A_{\ell'}=\E_{\ell}^{T}\A_{\ell}\E_{\ell}, \quad \E_{\ell}=\prod_{i=1}^{p_\ell} \E_{\I_{\ell,i}},
 \]
 where $\E_{\I_{\ell,i}}$ is defined following \eqref{V_schur}. The corresponding active DOFs are $\S_{\ell'} = \S_\ell \setminus \cup_{i = 1}^{p_\ell} \I_{\ell,i}$, which now comprise only edges and corners at this scale.
\item
 \textbf{Block Jacobi preconditioning.} Now group $\S_{\ell'}$ by edges and corners, and let $\{ \I_{\ell', i} \}_{i = 1}^{r_\ell}$ be the collection of corresponding index sets. Preconditioning as in Section \ref{sec:jacobi-precond} yields
 \[
  \A_{\ell''} = \C_{\ell}^T \A_{\ell'} \C_{\ell}, \quad \C_{\ell} = \prod_{i = 1}^{r_\ell} \C_{\I_{\ell', i}},
 \]
 where $\C_{\I_{\ell', i}}$ is as defined in \eqref{eqn:jacobi-precond}. The resulting matrix $\A_{\ell''}$ has unit block diagonal and the set of active DOFs $\S_{\ell''} = \S_{\ell'}$ is unchanged. This step is skipped (i.e., $\C_\ell = I$) in the standard HIF.
\item
 \textbf{Edge skeletonization.} Group $\S_{\ell''}$ by edges and let $\{ \I_{\ell'', i} \}_{i = 1}^{q_\ell}$ be the collection of corresponding index sets. Skeletonization as in Section \ref{skeletonization} then gives
 \[
  \A_{\ell+1}\approx\K_{\ell}^{T}\A_{\ell''}\K_{\ell}, \qquad
  \K_{\ell}=\prod_{i=1}^{q_\ell}\K_{\I_{\ell'',i}}, \quad
  \K_{\I_{\ell'',i}} =\Z_{\I_{\ell'',i}}\E_{\Ir_{\ell'',i}},
 \]
 where $\Z_{\I_{\ell'', i}}$ is defined following \eqref{eq_T} and $\E_{\Ir_{\ell'', i}}$ is the associated elimination matrix \eqref{V_schur} acting on the redundant indices. All DOFs $(\cup_{i = 1}^{p_\ell} \I_{\ell, i}) \cup (\cup_{i = 1}^{q_\ell} \Ir_{\ell'', i})$ have now been eliminated up to level $\ell$. The remaining active DOFs $\S_{\ell + 1} = \S_{\ell''} \setminus \cup_{i = 1}^{q_\ell} \Ir_{\ell'', i}$ consist of only edge skeletons and corners, with the former typically clustering around the latter.
\end{enumerate}

At the conclusion of this process, we have the final matrix
\[
 \A_L \approx R_{L-1}^T \cdots R_0^T \A R_0 \cdots R_{L-1}, \quad R_\ell = \E_\ell \C_\ell \K_\ell
\]
at the root, which is everywhere the identity except in the block indexed by the remaining active DOFs $\S_L$ containing, essentially, just the top-level corners. As a result, it is easily invertible, as is $R_\ell$ since each constituent factor is either block diagonal or triangular. Thus, we obtain the approximation $\F$ defined as follows
\begin{equation}
 \A \approx \F \equiv R_{0}^{-T} \cdots R_{L-1}^{-T}\A_{L}R_{L-1}^{-1} \cdots R_{0}^{-1},
 \label{eqA}
\end{equation}
and, by applying the inverse on both sides,
\begin{equation} 
 \A^{-1} \approx \F^{-1} = R_{0} \cdots R_{L-1}\A_{L}^{-1}R_{L-1}^{T} \cdots R_{0}^{T}.
 \label{eqAinv}
\end{equation}
The factorization
\begin{equation}
 \F = \G\G^T, \quad \G = R_0^{-T} \cdots R_{L-1}^{-T}B_L, \quad \A_L = B_L B_L^T
 \label{eqF}
\end{equation}
is a \defn{generalized Cholesky decomposition} (assuming that $\A_L$ is SPD), composed of a multilevel sequence of block sparse local matrices.

\begin{figure}
\centering \captionsetup[subfigure]{labelformat=empty}
\begin{subfigure}{0.225\textwidth}
 \centering
 \includegraphics[width=\textwidth,trim=1.65cm 1.65cm 0.75cm 1.15cm,clip]{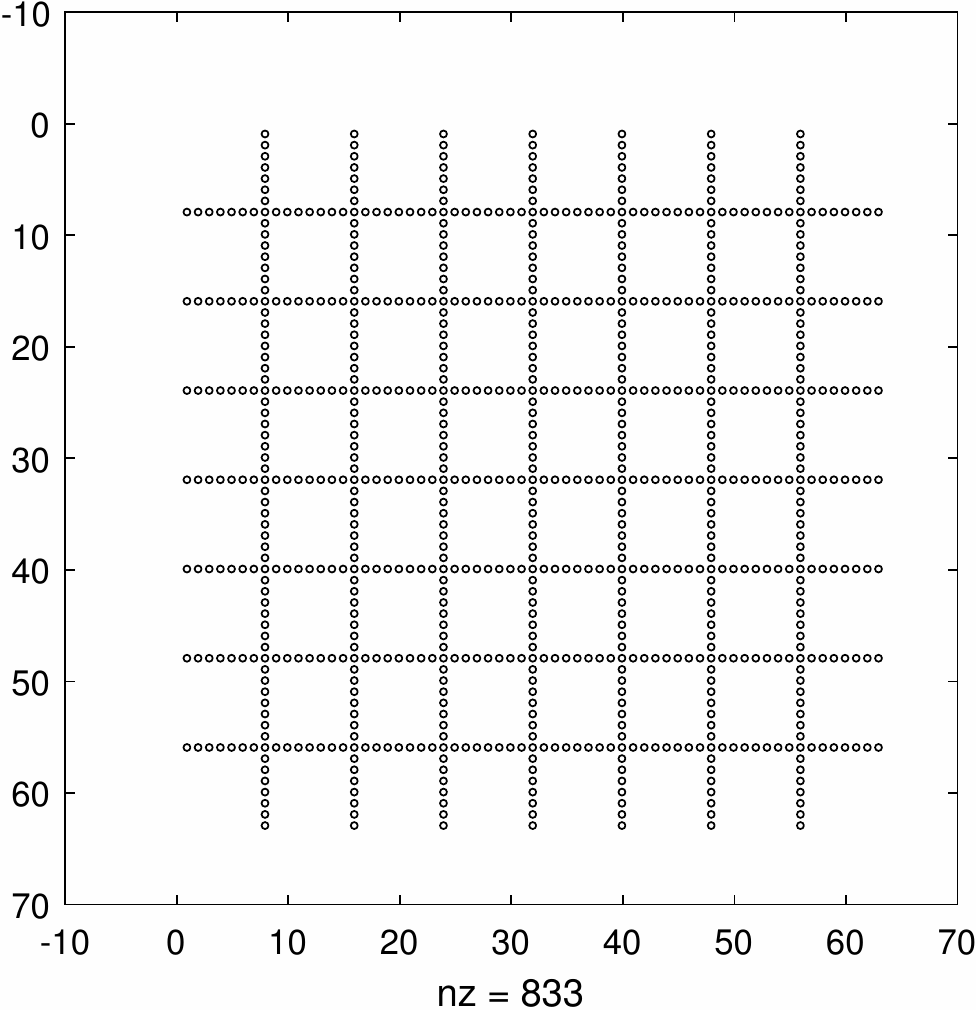}
 \caption{$\ell=0$ cell elimination}
\end{subfigure}
\begin{subfigure}{0.225\textwidth}
 \centering
 \includegraphics[width=\textwidth,trim=1.65cm 1.65cm 0.75cm 1.15cm,clip]{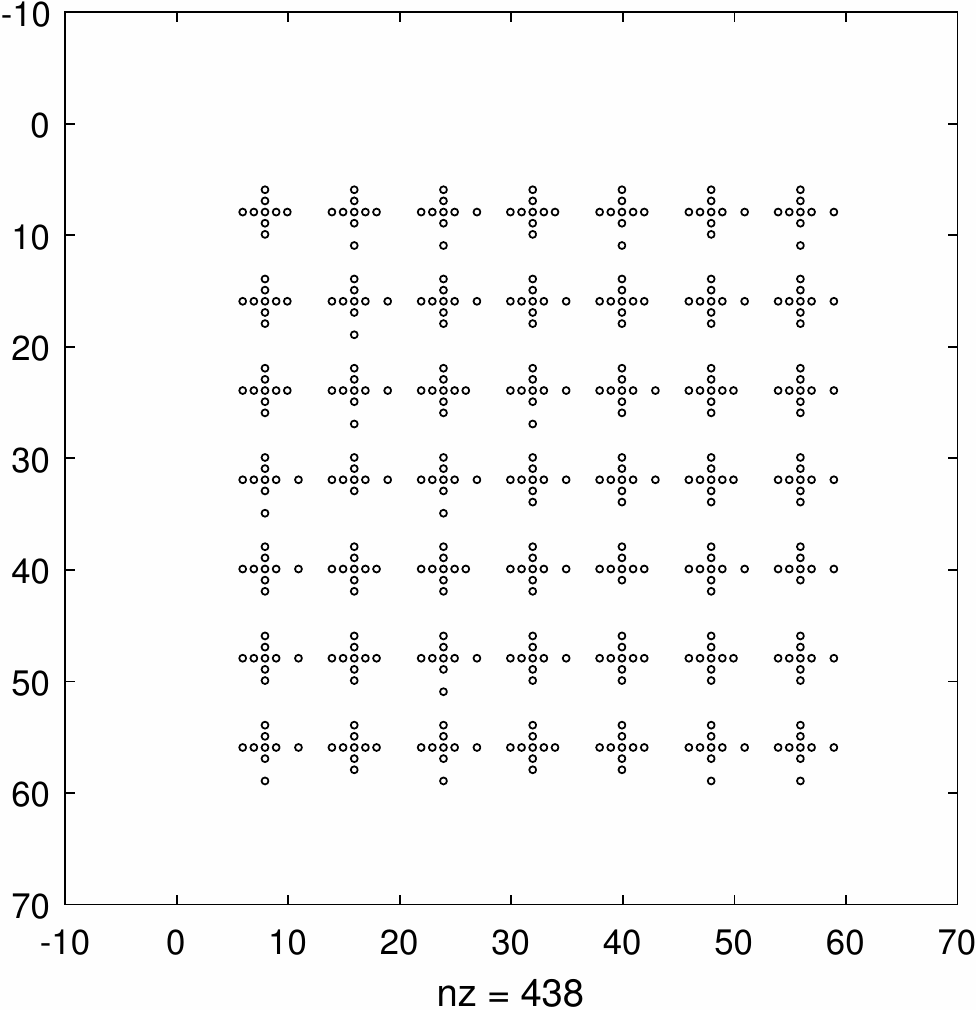}
 \caption{$\ell=0$ skeletonization}
\end{subfigure}
\begin{subfigure}{0.225\textwidth}
 \centering
 \includegraphics[width=\textwidth,trim=1.65cm 1.65cm 0.75cm 1.15cm,clip]{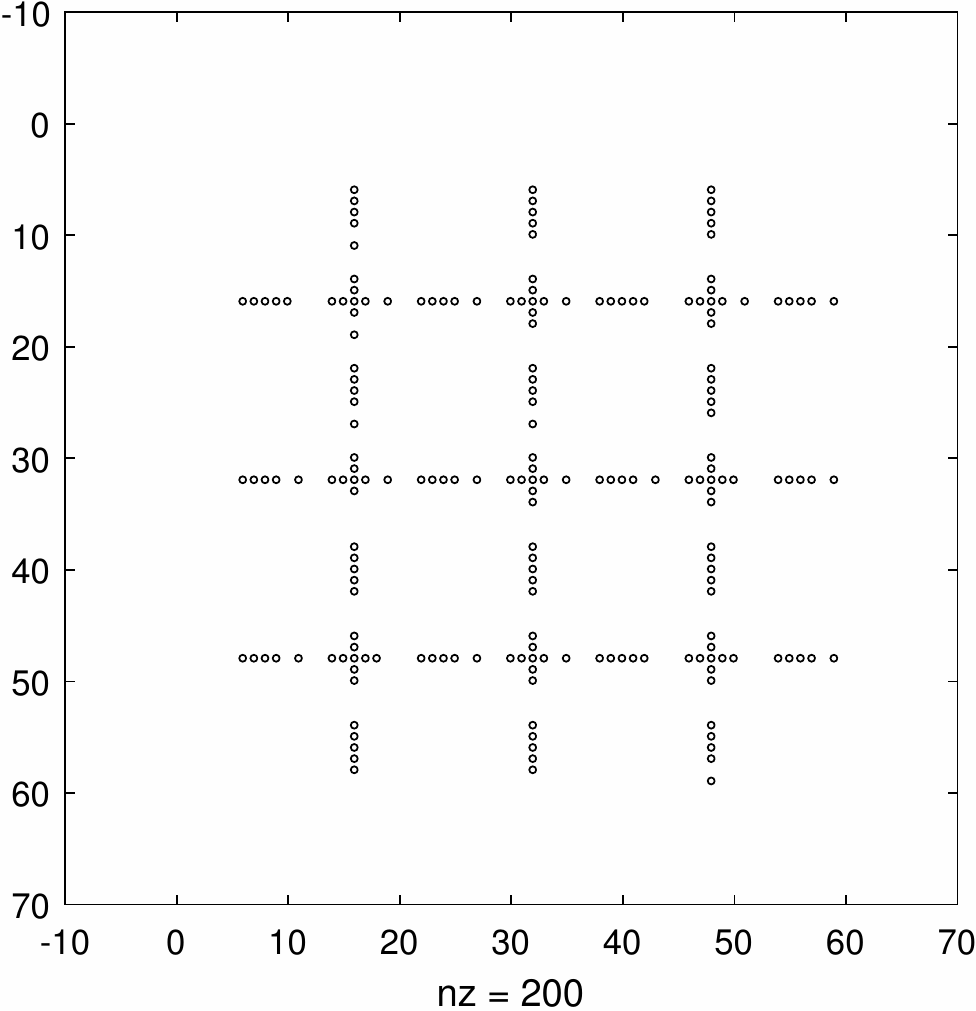}
 \caption{$\ell=1$ cell elimination}
\end{subfigure}\\
\begin{subfigure}{0.225\textwidth}
 \centering
 \includegraphics[width=\textwidth,trim=1.65cm 1.65cm 0.75cm 1.15cm,clip]{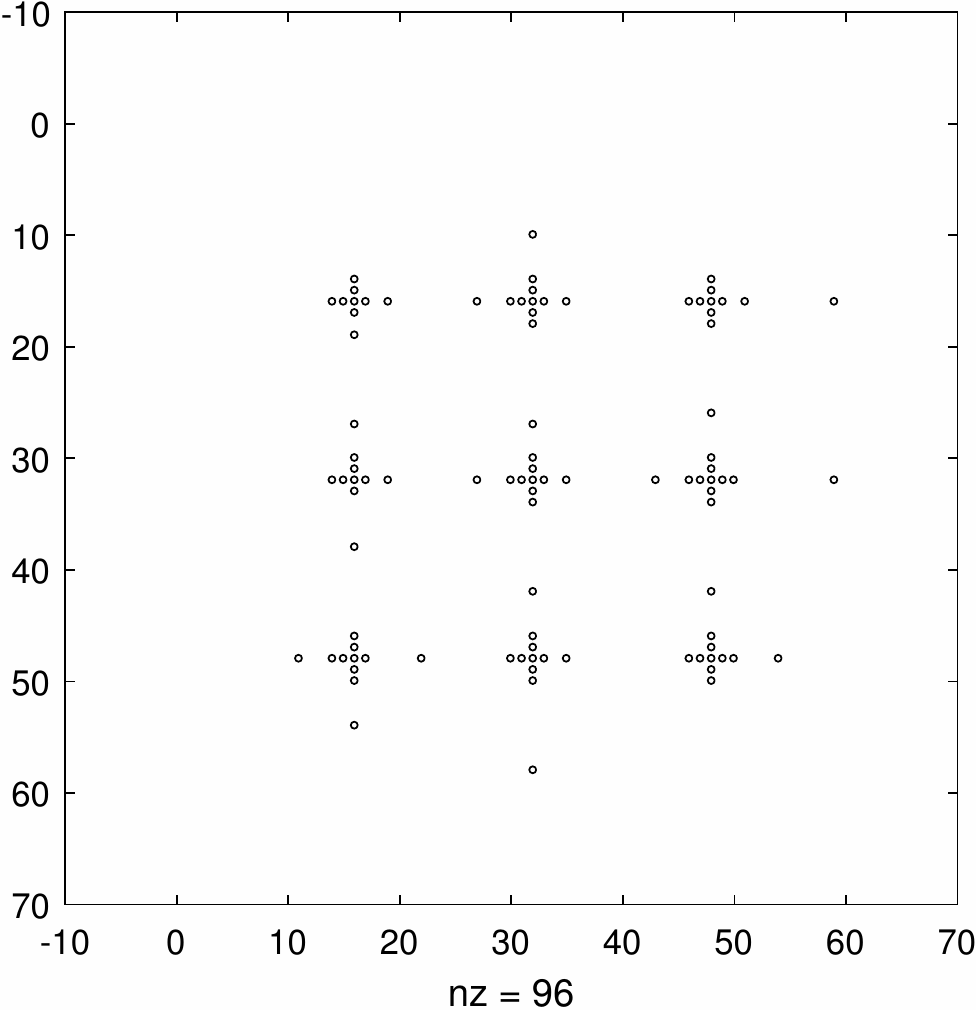}
 \caption{$\ell=1$ skeletonization}
\end{subfigure}
\begin{subfigure}{0.225\textwidth}
 \centering
 \includegraphics[width=\textwidth,trim=1.65cm 1.65cm 0.75cm 1.15cm,clip]{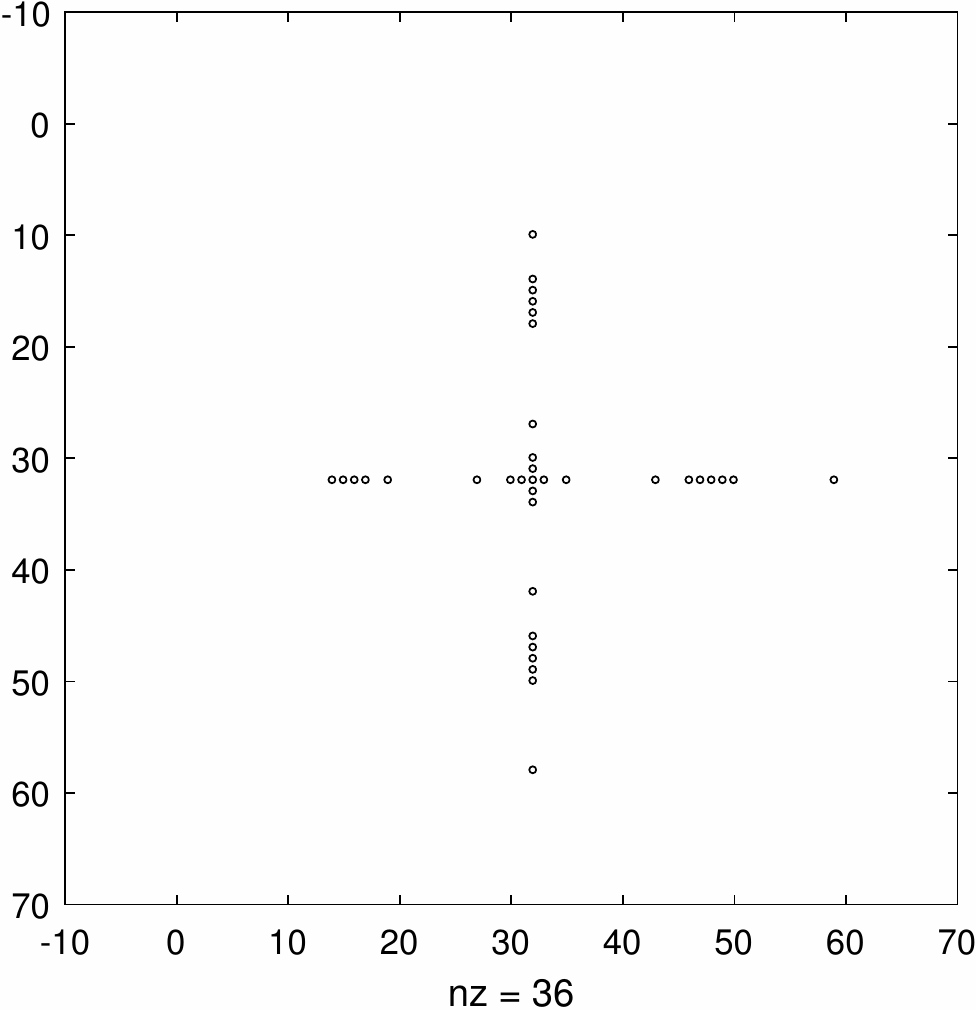}
 \caption{$\ell=2$ cell elimination}
\end{subfigure}
\begin{subfigure}{0.225\textwidth}
 \centering
 \includegraphics[width=\textwidth,trim=1.65cm 1.65cm 0.75cm 1.15cm,clip]{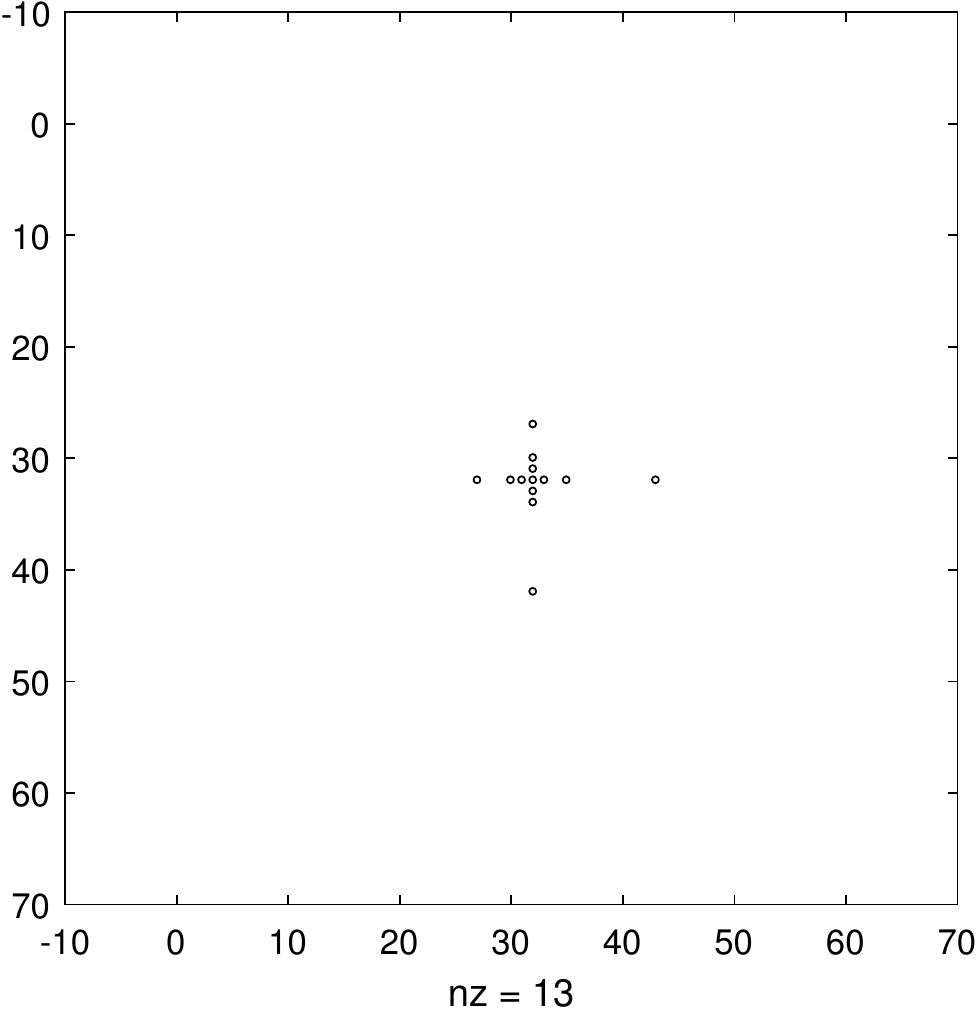}
 \caption{$\ell=2$ skeletonization}
\end{subfigure}
\caption{Active DOFs at each level $\ell$ of PHIF in 2D.}
\label{mHIFDE-fig}
\end{figure}

See Figure \ref{mHIFDE-fig} for an example of the elimination process.
Compared to HIF on the same problem (Figure \ref{HIFDE-fig}), the skeletons are somewhat more numerous and disperse, but generally retain the same structure. Furthermore, we observe that the skeletonization rank $\rho_\ell$ at each level $\ell$  still scales as $\rho_\ell = O(\ell)$ (see Section \ref{Numerical}) just as in HIF; consequently, the overall complexity is unchanged at $O(N)$ \cite[Theorem 4.6]{HIFDE}.

\subsection{Three dimensions}
\begin{figure}
\centering
\captionsetup[subfigure]{labelformat=empty}
\begin{subfigure}{0.225\textwidth}
 \centering
 \includegraphics[width=\textwidth,trim=2.1cm 2.1cm 2.1cm 1.9cm,clip]{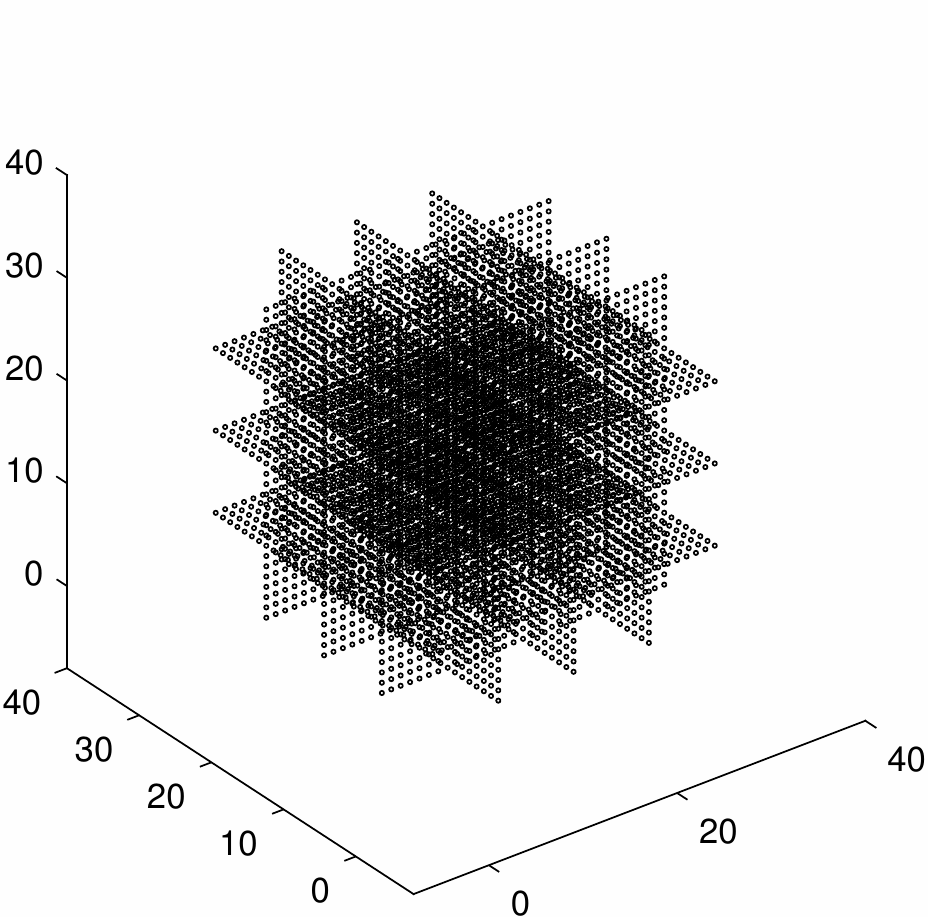}
 \caption{$\ell=0$ cell elimination}
\end{subfigure}
\begin{subfigure}{0.225\textwidth}
\centering
\includegraphics[width=\textwidth,trim=2.1cm 2.1cm 2.1cm 1.9cm,clip]{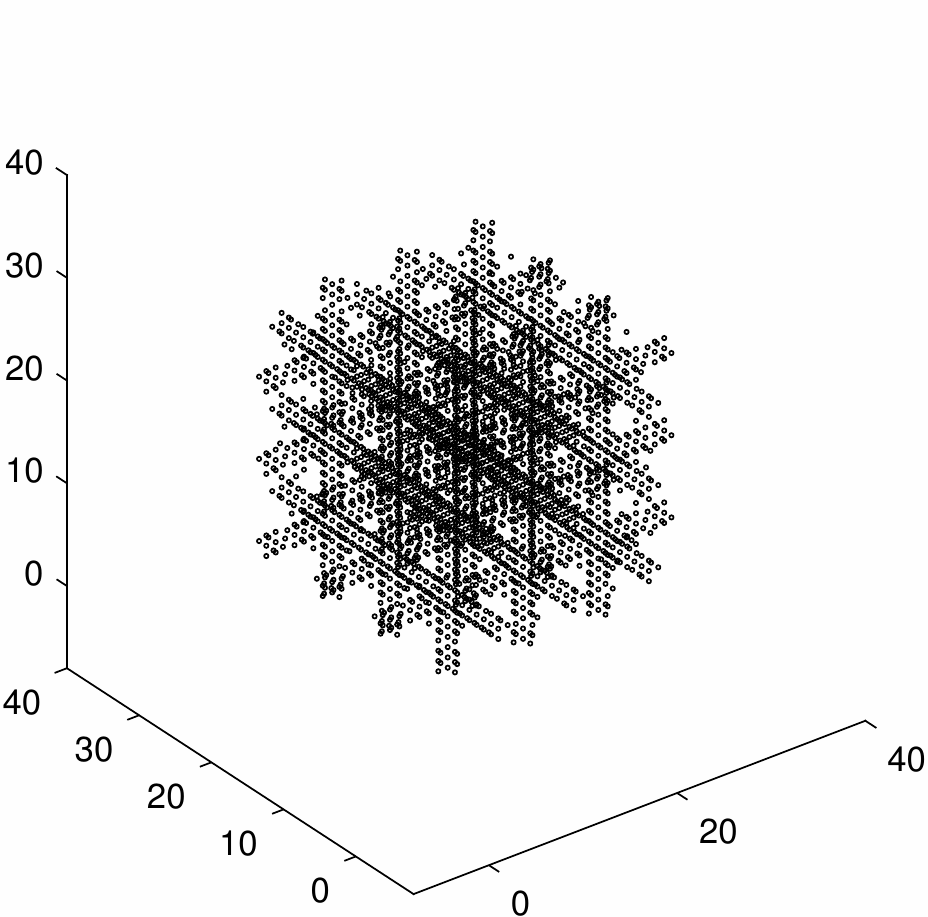}
\caption{$\ell=0$ skeletonization}
\end{subfigure}
\begin{subfigure}{0.225\textwidth}
 \centering
 \includegraphics[width=\textwidth,trim=2.1cm 2.1cm 2.1cm 1.9cm,clip]{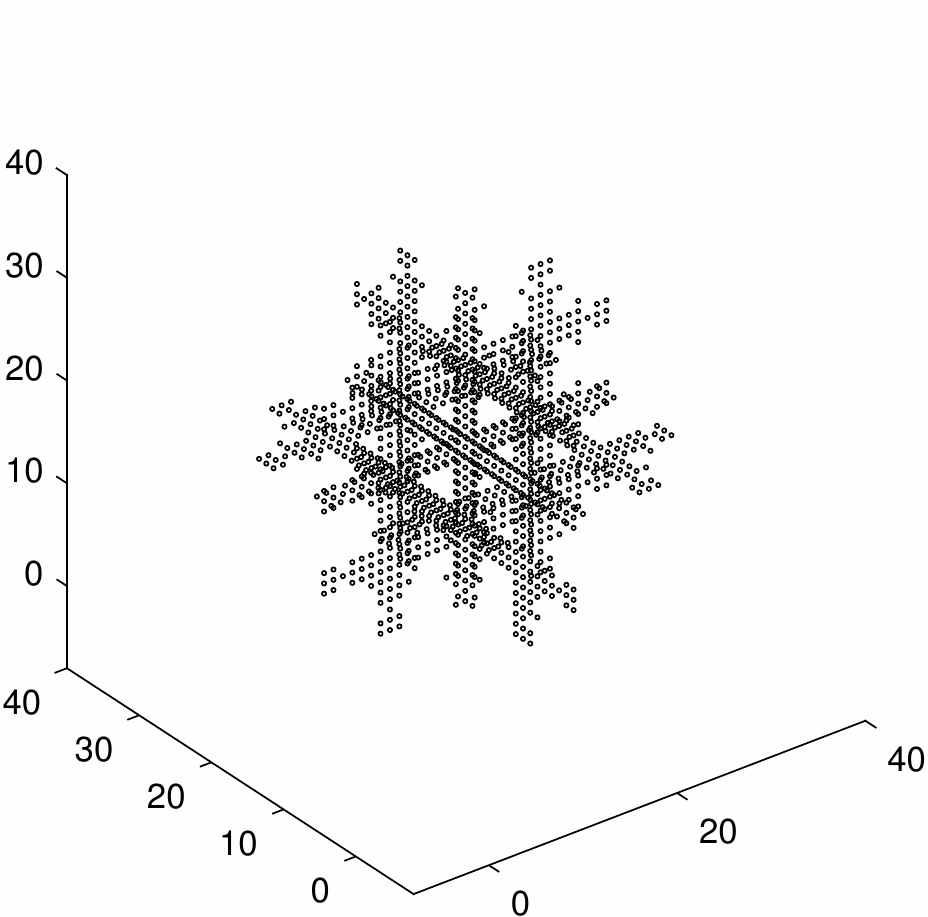}
 \caption{$\ell=1$ cell elimination}
\end{subfigure}
\begin{subfigure}{0.225\textwidth}
 \centering
 \includegraphics[width=\textwidth,trim=2.1cm 2.1cm 2.1cm 1.9cm,clip]{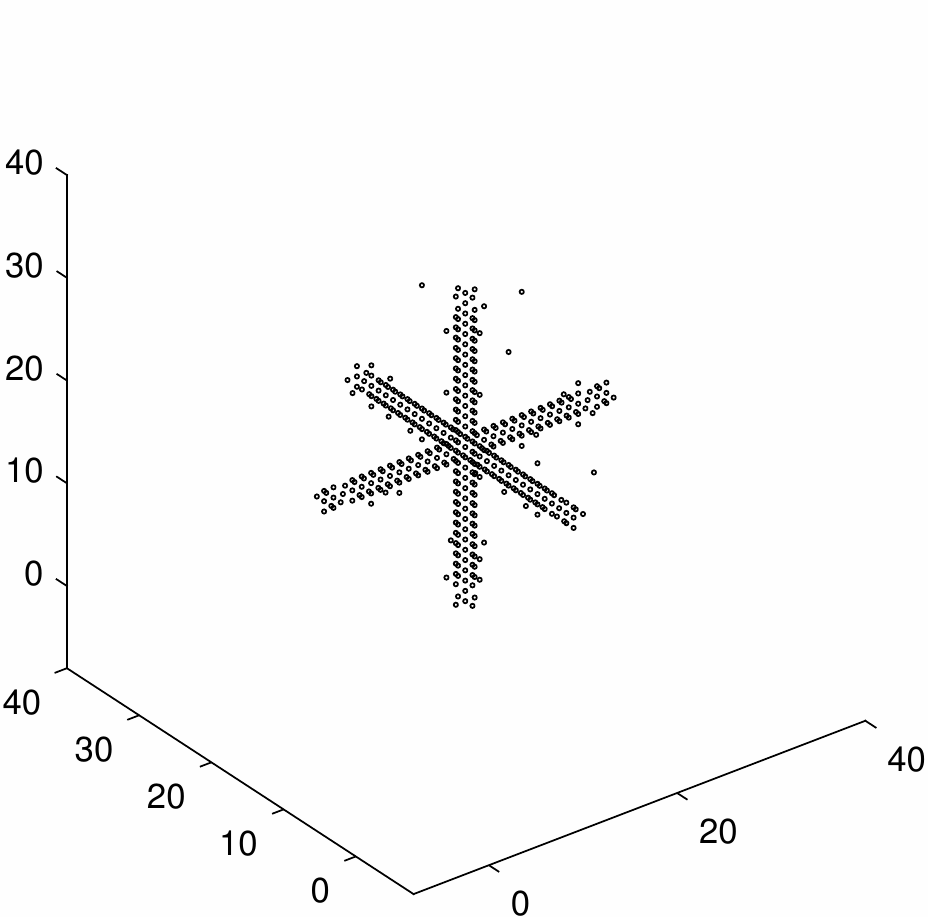}
 \caption{$\ell=1$ skeletonization}
\end{subfigure}
\caption{Active DOFs at each level $\ell$ of PHIF in 3D.}
\label{mHIFDE-fig3d}
\end{figure}

The algorithm can be extended to 3D in the natural way. Consider the analogous problem defined on $\Omega = (0,1)^3$ and discretized with the seven-point stencil. We now use an octree to hierarchically split $\Omega$ into $2^{3(L-\ell)}$ cubic cells at each level $\ell$, each covering $(2^{\ell} m + 1)^3$ points. Among these,
\begin{itemize}
\item $(2^{\ell} m - 1)^3$ are \defn{interior} points;
\item $6(2^{\ell} m - 1)$ are \defn{face} points, where each face is shared between up to two cells;
\item $12(2^{\ell} m - 1)$ are \defn{edge} points, where each edge is shared between up to four cells; and
\item $8$ are \defn{corner} points, each shared between up to eight cells.
\end{itemize}
Denote by $p_\ell$ the total number of cells on level $\ell$, $q_\ell \sim 3 p_\ell$ the total number of faces, and $r_\ell \sim 7 p_\ell$ the total number of faces, edges, and corners altogether. Then we analogously perform at each level $\ell$:
\begin{enumerate}
\item interior cell elimination for each of $p_\ell$ DOF sets;
\item face, edge, and corner preconditioning for each of $r_\ell$ DOF sets; and
\item face skeletonization for each of $q_\ell$ DOF sets.
\end{enumerate}
The remaining active DOFs tend to cluster around edges and corners. See Figure \ref{mHIFDE-fig3d} for an example; as in 2D, the skeletons are slightly less ordered than for HIF (Figure \ref{HIFDE-fig3d}).
The skeletonization rank is $\rho_\ell = O(2^\ell)$ just as before, for a total cost of $O(N \log N)$ \cite{HIFDE}.

\begin{remark}
\label{rmk:3d-edge-skel}
In HIF, estimated $O(N)$ complexity can be achieved by additional edge skeletonization after face skeletonization to reduce the rank back down to $\rho_\ell = O(\ell)$. The same strategy presumably translates to PHIF, where we now would employ two extra steps:
\begin{enumerate}
\item[4.] edge and corner preconditioning, and
\item[5.] edge skeletonization.
\end{enumerate}
Note that we did not explore this approach here.
\end{remark}

\begin{remark}
As presently formulated, the ID-based skeletonization in Section \ref{skeletonization} is actually not critical and other low-rank compression techniques may well be used, including the SVD \cite{HIFDE}. The ID is required only for the edge skeletonization variant of Remark \ref{rmk:3d-edge-skel}, which depends on the geometry being preserved.
\end{remark}

\subsection{Heuristic error analysis}\label{sec:error}
In order to understand the impact of PHIF on the solve error, let us write each level of the algorithm more carefully as $\A_{\ell + 1} = R_\ell^T \A_\ell R_\ell + E_\ell$, where $\| E_\ell \| = O(\epsilon \| \A_\ell \|)$ is the approximation error from skeletonization to some relative precision $\epsilon$. To ease the notation, we define a chain of matrix products with an ellipsis, e.g., $R_{L-1}^T \cdots R_{\ell + 1}^T$, where the whole chain of products is taken to be the identity if the subscript of some matrix in the chain is greater than $L-1$. If both matrices in the endpoints are the same, the chain reduces to just one matrix, e.g., $R_{L-1}^T \cdots R_{L-1}^T = R_{L-1}^T$. Then 
\begin{align*}
\A_L &= R_{L-1}^T \A_{L-1} R_{L-1} + E_{L-1}\\
&= R_{L-1}^T (R_{L-2}^T \A_{L-2} R_{L-2} + E_{L-2}) R_{L-1} + E_{L-1}\\
&= R_{L-1}^T \cdots R_0^T \A_0 R_0 \cdots R_{L-1} + \sum_{\ell = 0}^{L-1} R_{L-1}^T \cdots R_{\ell + 1}^T E_\ell R_{\ell + 1} \cdots R_{L-1},
\end{align*}
so
\[
\A = \A_0 = R_0^{-T} \cdots R_{L-1}^{-T} \A_L R_{L-1}^{-1} \cdots R_0^{-1} - \sum_{\ell = 0}^{L-1} R_0^{-T} \cdots R_\ell^{-T} E_\ell R_\ell^{-1} \cdots R_0^{-1},
\]
with approximate factorization
\[
\F = \G \G^T, \quad \G = R_0^{-T} \cdots R_{L-1}^{-T} B_L, \quad \A_L = B_L B_L^T.
\]
Hence the solve error is 
\begin{subequations}
\label{eqn:solve-err1}
\begin{align}
\| I - \G^{-1} \A \G^{-T} \| &= \| B_L^{-1} \| \left\| \sum_{\ell = 0}^{L-1} R_{L-1}^T \cdots R_{\ell + 1}^T E_\ell R_{\ell + 1} \cdots R_{L-1} \right\| \| B_L^{-T} \|\\
&\lesssim \epsilon \| \A_L^{-1} \| \sum_{\ell = 0}^{L-1} \| \A_\ell \| \| R_{\ell + 1} \cdots R_{L-1} \|^2.
\end{align}
\end{subequations}
To bound $\| A_\ell \|$, note that
\[
A_\ell = R_\ell^{-T} \cdots R_{L -1}^{-T} A_L R_{L - 1}^{-1} \cdots R_\ell^{-1} + O(\epsilon),
\]
so
\[
\| A_\ell \| \lesssim \| A_L \| \| R_\ell^{-1} \|^2 \| (R_{\ell + 1} \cdots R_{L - 1})^{-1} \|^2.
\]
Combining with \eqref{eqn:solve-err1} yields
\begin{equation}
\| I - G^{-1} A G^{-T} \| \lesssim \epsilon \kappa (A_L) \sum_{\ell = 0}^{L - 1} \left[ \| R_\ell^{-1} \| \: \kappa (R_{\ell + 1} \cdots R_{L - 1}) \right]^2.
\label{eqn:solve-err2}
\end{equation}

The estimate \eqref{eqn:solve-err2} holds for both HIF and PHIF, given appropriate interpretation of $R_\ell$ and $A_\ell$. We now contrast the two approaches by estimating and directly comparing each term.

In HIF, $R_\ell = \E_\ell \K_\ell$, where:
\begin{itemize}
\item
$\E_{\ell}=\prod_{i=1}^{p_\ell} \E_{\I_{\ell,i}}$ with, referring to Section \ref{preliminaries}, the norm of each matrix and its inverse bounded as
\[
\|\E_{\I_{\ell,i}}^{\pm 1}\| \leq \|L_{\I_{\ell,i}}^{\mp 1}\| \left( 1 + \|L_{\I_{\ell,i}}^{-1} (\A_\ell)_{\B_{\ell,i}\I_{\ell,i}}^{T}\| \right).
\]
Since $\E_\ell$ is nearly block diagonal, with coupling between the $\E_{\I_{\ell,i}}$ only along shared interfaces, we can expect
\begin{subequations}
\label{eqn:cond-elim}
\begin{align}
\kappa (\E_\ell) &\lesssim \max_{i, j} \|\L_{\I_{\ell,i}}\| \|\L_{\I_{\ell,j}}^{-1}\| \left( 1 + \|L_{\I_{\ell,i}}^{-1} (\A_\ell)_{\B_{\ell,i}\I_{\ell,i}}^{T}\| \right) \left( 1 + \|L_{\I_{\ell,j}}^{-1} (\A_\ell)_{\B_{\ell,j}\I_{\ell,j}}^{T}\| \right)\\
&= \max_{i, j} \frac{\| \L_{\I_{\ell, i}} \|}{\| \L_{\I_{\ell, j}} \|} \: \kappa (\L_{\I_{\ell, j}}) \left( 1 + \|L_{\I_{\ell,i}}^{-1} (\A_\ell)_{\B_{\ell,i}\I_{\ell,i}}^{T}\| \right) \left( 1 + \|L_{\I_{\ell,j}}^{-1} (\A_\ell)_{\B_{\ell,j}\I_{\ell,j}}^{T}\| \right).
\end{align}
\end{subequations}
In other words, $\kappa (\E_\ell)$ depends on the uniformity as well as the conditioning of the $\L_{\I_{\ell, i}}$.
\item
Similarly, $\K_\ell = \prod_{i = 1}^{q_\ell} \K_{\I_{\ell'', i}}$ (following the above algorithm with $\C_\ell = I$), where $\K_{\I_{\ell'', i}} = \Z_{\I_{\ell'', i}} \E_{\Ir_{\ell'', i}}$. But $\| \Z_{\I_{\ell'', i}}^{\pm 1} \| \leq 1 + \| \T_{\I_{\ell'', i}} \|$, where, in practice, $\| \T_{\I_{\ell', i}} \|$ is small due to the stability of the ID \cite{ID}, so $\| \Z_{\I_{\ell'', i}}^{\pm 1} \| = O(1)$ (or growing slowly with the size of the matrix). Hence, $\| \K_{\I_{\ell'', i}}^{\pm 1} \| \lesssim \| \E_{\Ir_{\ell'', i}}^{\pm 1} \|$. Now use that $\K_\ell$ is block diagonal since each $\K_{\I_{\ell'', i}}$ acts on disjoint separators to obtain
\begin{subequations}
\label{eqn:cond-skel}
\begin{align}
\kappa (\K_\ell) &\lesssim \max_{i, j} \| \E_{\Ir_{\ell'', i}} \| \| \E_{\Ir_{\ell'', j}}^{-1} \|\\
&\leq \frac{\| \L_{\Ir_{\ell'', i}} \|}{\| \L_{\Ir_{\ell'', j}} \|} \: \kappa (\L_{\Ir_{\ell'', j}}) \left( 1 + \| \L_{\Ir_{\ell'', i}}^{-1} (\widetilde{\A_{\ell''}})_{\Is_{\ell'', i} \Ir_{\ell'', i}}^T \| \right) \left( 1 + \| \L_{\Ir_{\ell'', j}}^{-1} (\widetilde{\A_{\ell''}})_{\Is_{\ell'', j} \Ir_{\ell'', j}}^T \| \right).
\end{align}
\end{subequations}
\end{itemize}

On the other hand, in PHIF, $R_\ell = \E_\ell \C_\ell \K_\ell$, where:
\begin{itemize}
\item
$\E_0$ at the initial level is the same as in HIF. At all subsequent levels, however, $\E_\ell$ acts on an $\A_\ell$ that has been preconditioned by $\C_{\ell'}$ for all $\ell' < \ell$. This has the effect of driving $\| \L_{\I_{\ell, i}} \|$ towards $1$ (since the corresponding $(\A_\ell)_{\I_{\ell, i} \I_{\ell, i}}$ has been rescaled to $O(1)$) as well as reducing $\kappa (\L_{\I_{\ell, i}})$. Furthermore, if $\A$ is SPD, then the off-diagonal term $\| (\A_\ell)_{\B_{\ell, i} \I_{\ell, i}} \|$ also cannot be too large. Thus, we can expect each quantity in \eqref{eqn:cond-elim} to be better behaved in PHIF than in HIF.
\item
$\C_\ell = \prod_{i = 1}^{r_\ell} \C_{\ell, i}$, where clearly $\kappa (\C_\ell) = \max_{i,j} \| \L_{\I_{\ell', i}}\|\| \L_{\I_{\ell', j}}^{-1}\|$. As with $\E_\ell$, each $C_\ell$ for $\ell > 0$ is typically not too ill-conditioned due to preconditioning at previous levels. In effect, this is similar to the $\kappa (\K_\ell)$ term in HIF but can be much better for $\ell > 0$ due to prior preconditioning.
\item
Again, $\kappa (\K_\ell)$ can be bounded as in \eqref{eqn:cond-skel}, but now with each term well-behaved following the same remarks as for $\E_\ell$. In fact, we argue that $\kappa (\K_\ell) = O(1)$ since the post-preconditioned submatrix on which each elimination operator acts has the form
\[
\begin{bmatrix}
(\A_{\ell''})_{\Ir_{\ell'', i} \Ir_{\ell'', i}} & (\A_{\ell''})_{\Is_{\ell'', i} \Ir_{\ell'', i}}^T\\
(\A_{\ell''})_{\Is_{\ell'', i} \Ir_{\ell'', i}} & (\A_{\ell''})_{\Is_{\ell'', i} \Is_{\ell'', i}}
\end{bmatrix} =
\begin{bmatrix}
I + \T_{\I_{\ell'', i}}^T \T_{\I_{\ell'', i}} & -\T_{\I_{\ell'', i}}^T\\
-\T_{\I_{\ell'', i}} & I
\end{bmatrix},
\]
where all terms are $O(1)$.
\end{itemize}
Putting these together, we have that
\[
\frac{\kappa' (R_\ell)}{\kappa (R_\ell)} \sim \frac{\kappa' (\E_\ell)}{\kappa (\E_\ell)} \: \kappa' (\C_\ell) \: \frac{\kappa' (\K_\ell)}{\kappa (\K_\ell)} \lesssim \frac{\kappa' (\C_\ell)}{\kappa (\K_\ell)} \lesssim 1,
\]
where $\kappa (\cdot)$ and $\kappa' (\cdot)$ denote, respectively, quantities for HIF and PHIF. Looking at \eqref{eqn:solve-err2}, it remains to compare $\| R_\ell^{-1} \|$ and $\kappa (\A_L)$. For the former, a simplified analysis using $\| \L_{\I_{\ell, i}} \|\sim O(1)$ in PHIF and omitting some lengthy terms of the form $(1+\|L^{-1}A^T\|)$, which are worse in HIF than PHIF, gives
\[
\| R_\ell^{-1} \| \lesssim
\begin{cases}
 \max_{i,j} \| \L_{\I_{\ell, i}} \| \| \L_{\Ir_{\ell'', j}} \| & \text{(HIF)}\\
\max_i \| \L_{\I_{\ell', i}} \| & \text{(PHIF)}
\end{cases}
\]
which we argue is smaller for PHIF than HIF by similar reasoning as for comparing $\C_\ell$ in PHIF with $\K_\ell$ in HIF above. Finally, we have $\kappa (\A_L)$, which, of course, is typically much improved in PHIF since it has been preconditioned throughout.

This analysis provides some insight on the essential feature of the algorithm: whereas $\A_L$ can be very ill-conditioned in HIF, inheriting in large part from $\A$ itself, in PHIF it results from multiple levels of preconditioning through the action of the rescaling operators $\C_\ell$ and so should be better behaved---in certain cases, much better, by several orders of magnitude, as we will see in Section \ref{Numerical}.

It is also instructive to revisit Remark \ref{rmk:skeleton-spd} on the required accuracy to remain SPD in light of PHIF. Applying the remark, the intermediate matrices $\A_\ell$ are all SPD and so the algorithm succeeds if $\| E_\ell \| < \lambda_{\min} (\A_\ell) = 1 / \| \A_\ell^{-1} \|$; for $\| E_\ell \| = O(\epsilon \| \A_\ell \|)$ as above, this gives $\epsilon \lesssim 1 / \max_\ell \kappa (\A_\ell)$. Naturally, $\kappa (\A_\ell)$ can be much smaller in PHIF than HIF, thereby significantly alleviating this issue.

\section{Numerical results}\label{Numerical}
We now present some benchmark examples to demonstrate the performance of PHIF on strongly ill-conditioned problems. Specifically, we consider the PDE \eqref{eq1} on $\Omega = (0, 1)^d$, discretized as in Section \ref{Algorithm} with $b(x) \equiv 0$ and $a(x)$ a quantized high-contrast random field defined as follows:
\begin{enumerate}
\item
 Initialize by sampling each grid point $a_{j} = a(x_j)$ from the standard uniform distribution.
\item
 Convolve with an isotropic Gaussian of width $4h$ to create some correlation structure.
\item
 Quantize by setting
 \[
  a_{j} =
  \begin{cases}
   10^{-2}, & a_{j}\le\mu,\\
   10^{2},  & a_{j}>\mu,
  \end{cases}
  \]
  where $\mu$ is the median of $\{a_{j}\}$.
\end{enumerate}
The resulting matrix $\A$ in \eqref{eq2} has condition number $\kappa (\A) = O(\sigma h^{-2}) = O(\sigma N^{2/d})$, where $\sigma = 10^4$ is the contrast ratio.

In the following, we report
\begin{itemize}
\item $\epsilon$: relative precision for ID compression;
\item $N$: total number of DOFs in the problem;
\item $|\S_L|$: number of active DOFs remaining at the highest level;
\item $e_a$: estimated apply error $\|\A-\F\|/\|\A\|$;
\item $e_s$: estimated solve error $\| I - \G^{-1} \A \G^{-T} \| \sim \| I - \A \F^{-1} \| \geq \| \A^{-1} - \F^{-1} \| / \| \A^{-1} \|$; and
\item $n_i$: number of CG iterations using $\F^{-1}$ as a preconditioner to solve to a relative residual of $10^{-12}$.
\end{itemize}
The errors $e_a$ and $e_s$ are estimated using randomized power iteration \cite{Dixon,Kuczynski} to $10^{-2}$ relative precision. All numerical experiments are performed in MATLAB using codes modified from FLAM \cite{FLAM}.

\subsection{Two dimensions}
\begin{figure}
\centering
\begin{subfigure}{0.4\textwidth}
 \includegraphics[width=\textwidth]{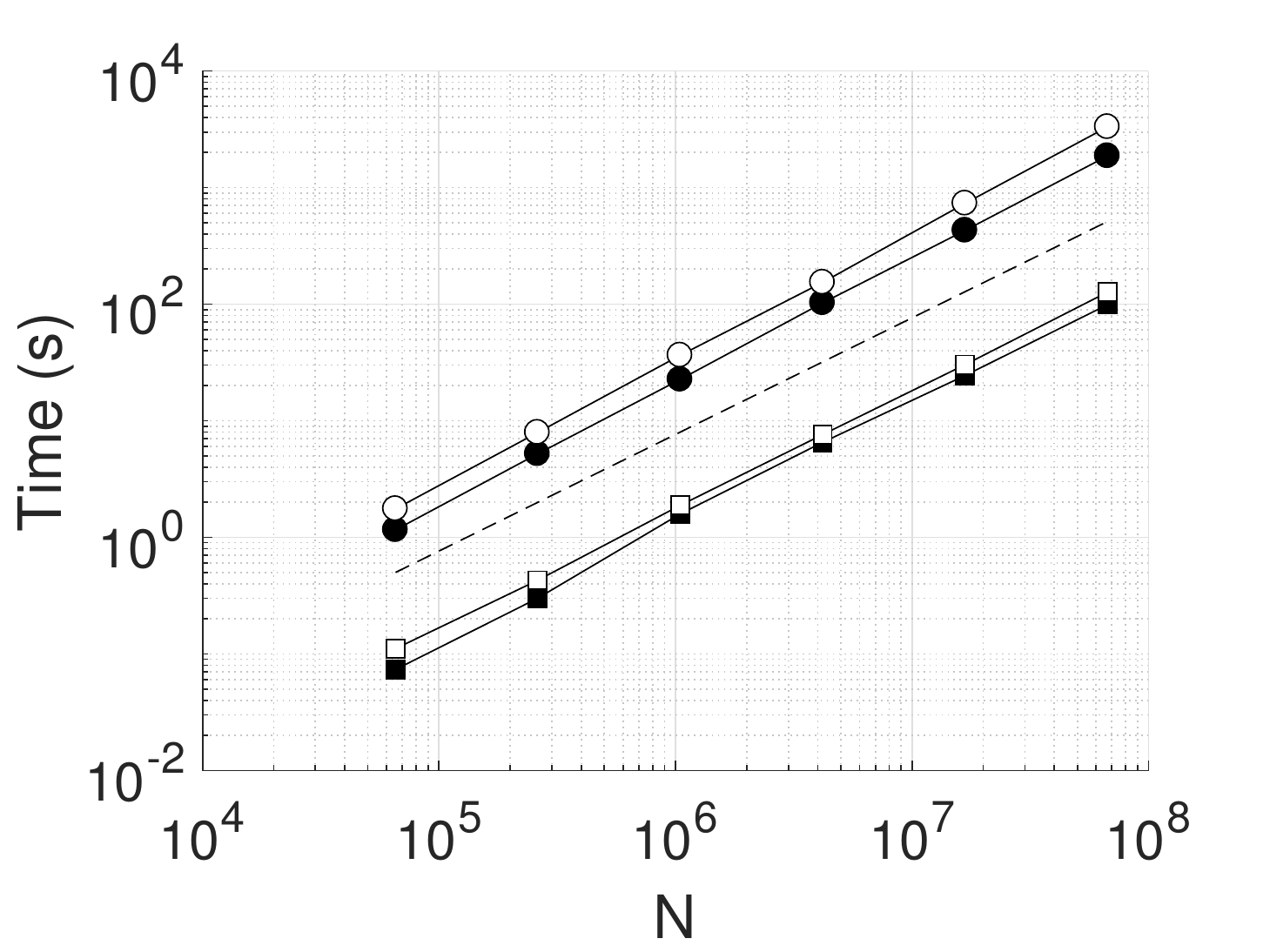}
\end{subfigure}
\hspace{2em}
\begin{subfigure}{0.4\textwidth}
\includegraphics[width=\textwidth]{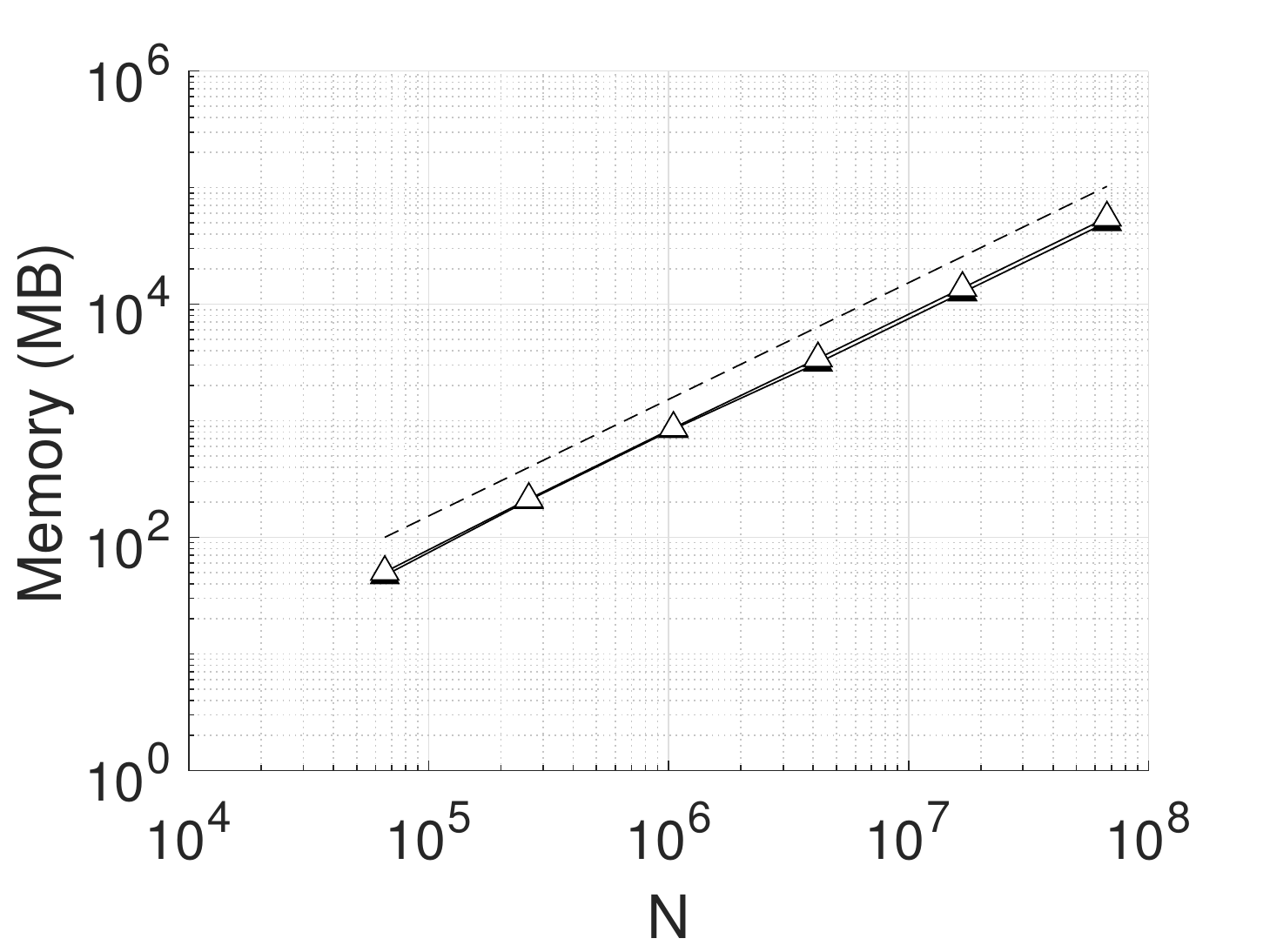}
\end{subfigure}
\caption{Scaling results for 2D example, showing factorization times ($\circ)$, application times for $\F^{-1}b$ ($\Box$), and memory storage ($\bigtriangleup$) for HIF (black) and PHIF (white) at $\epsilon=10^{-6}$. Reference $O(N)$ lines (dashed) are also included.}
\label{2dcomplexity}
\end{figure}

Consider first the example in 2D. Numerical results are given in Table \ref{t_results} with scaling plots shown in Figure \ref{2dcomplexity}. We immediately see the effectiveness of PHIF at improving $e_s$ (by $10^2$--$10^3$) while maintaining comparable $e_a = O(\epsilon)$ with HIF. This directly manifests in a smaller $n_i$ across all cases tested. We also observe the intermediate matrices encountered throughout PHIF to be much better conditioned, in agreement with Section \ref{sec:error}. Likewise, PHIF exhibits greater robustness with respect to remaining SPD; indeed, while HIF fails for $\epsilon = 10^{-4}$, PHIF is still able to produce good preconditioners at the same tolerance.

\begin{table}[htb]
\centering
\caption{Numerical results for 2D example. HIF at $\epsilon = 10^{-4}$ fails due to loss of positive definiteness.}
\label{t_results}
\begin{tabular}{cc|cccc|cccc}
\toprule
&& \multicolumn{4}{|c|}{HIF} & \multicolumn{4}{|c}{PHIF}\\
$\epsilon$ & $N$ & $|\S_L|$ & $e_a$ & $e_s$ & $n_i$ & $|\S_L|$ & $e_a$ & $e_s$ & $n_i$\\
\midrule
\multirow{4}{*}{$10^{-4}$} & $1023^2$ & \multicolumn{4}{|c|}{\multirow{4}{*}{not SPD}} & $60$ & $4.7$e$-5$ & $1.4$e$-1$ & $9$\\
& $2047^2$ &&&&& $68$ & $6.0$e$-5$ & $1.6$e$-1$ & $12$\\
& $4095^2$ &&&&& $77$ & $7.6$e$-5$ & $2.4$e$-1$ & $14$\\
& $8191^2$ &&&&& $86$ & $7.7$e$-5$ & $5.7$e$-1$ & $17$\\
\midrule
\multirow{4}{*}{$10^{-6}$} & $1023^2$ & $59$ & $2.9$e$-6$ & $7.3$e$-1$ & $16$ & $81$ & $4.9$e$-7$ & $1.1$e$-3$ & $4$\\
& $2047^2$ & $61$ & $2.9$e$-6$ & $8.5$e$-1$ & $20$ & $91$ & $6.9$e$-7$ & $1.5$e$-3$& $4$\\
& $4095^2$ & $63$ & $4.2$e$-6$ & $9.3$e$-1$ & $32$ & $111$ & $8.5$e$-7$ & $2.2$e$-3$ & $5$\\
& $8191^2$ & $46$ & $5.9$e$-6$ & $9.7$e$-1$ & $54$ & $125$ & $9.8$e$-7$ & $3.5$e$-3$ & $5$\\
\midrule
\multirow{4}{*}{$10^{-8}$} & $1023^2$ & $80$ & $3.1$e$-8$ & $2.9$e$-3$ & $4$ & $105$ & $6.5$e$-9$ & $7.1$e$-6$ & $4$\\
& $2047^2$ & $82$ & $3.2$e$-8$ & $5.3$e$-3$ & $5$ & $118$ & $8.2$e$-9$ & $1.8$e$-5$ & $3$\\
& $4095^2$ & $103$ & $3.5$e$-8$ & $1.1$e$-2$ & $5$ & $133$ & $9.7$e$-9$ & $3.1$e$-5$ & $3$\\
& $8191^2$ & $112$ & $4.4$e$-8$ & $2.1$e$-2$ & $6$ & $154$ & $1.1$e$-8$ & $2.8$e$-5$ & $3$\\
\bottomrule
\end{tabular}
\end{table}

As with HIF, PHIF achieves linear complexity, as is evident from the figure and from verifying the mild growth of $|\S_L|$ with $N$. However, at fixed $\epsilon$, PHIF can be quite a bit more expensive, due primarily to the extra Cholesky preconditioning steps. For example, at $N = 8191^2$ and $\epsilon = 10^{-6}$, these (1456 s) account for about 44\% of the total PHIF factorization time (3317 s) or 78\% of the equivalent HIF time (1869 s); in other words, just the multilevel preconditioning itself can almost double the factorization cost without even considering the potential added impact of larger skeleton sizes. Despite this overhead, PHIF is still able to decrease the total time to solution: for the same case, after factorization, HIF takes 5698 s to run $n_i = 54$ iterations for a total of 7567 s, while PHIF needs only 1199 s over $n_i = 5$ for a total of 4516 s --- a reduction of close to 40\%.

In a sense, the behavior of PHIF is similar to that of HIF computed at a higher precision: $e_s$ is improved at the cost of increased ranks. But careful inspection of the results suggests that the interpretation is not quite so simple. Comparing, e.g., $N = 1023^2$ for PHIF at $\epsilon = 10^{-4}$ vs.\ HIF at $\epsilon = 10^{-6}$ and for PHIF at $\epsilon = 10^{-6}$ vs.\ HIF at $\epsilon = 10^{-8}$ reveals that PHIF can use approximately the same final rank to achieve better $e_s$. At the same time, $e_a$ is slightly worse, thus indicating a certain balance between the forward and inverse errors. We will see such nuances more clearly in the next example as well.

\subsection{Three dimensions}
\begin{figure}
\centering
\begin{subfigure}{0.4\textwidth}
\includegraphics[width=\textwidth]{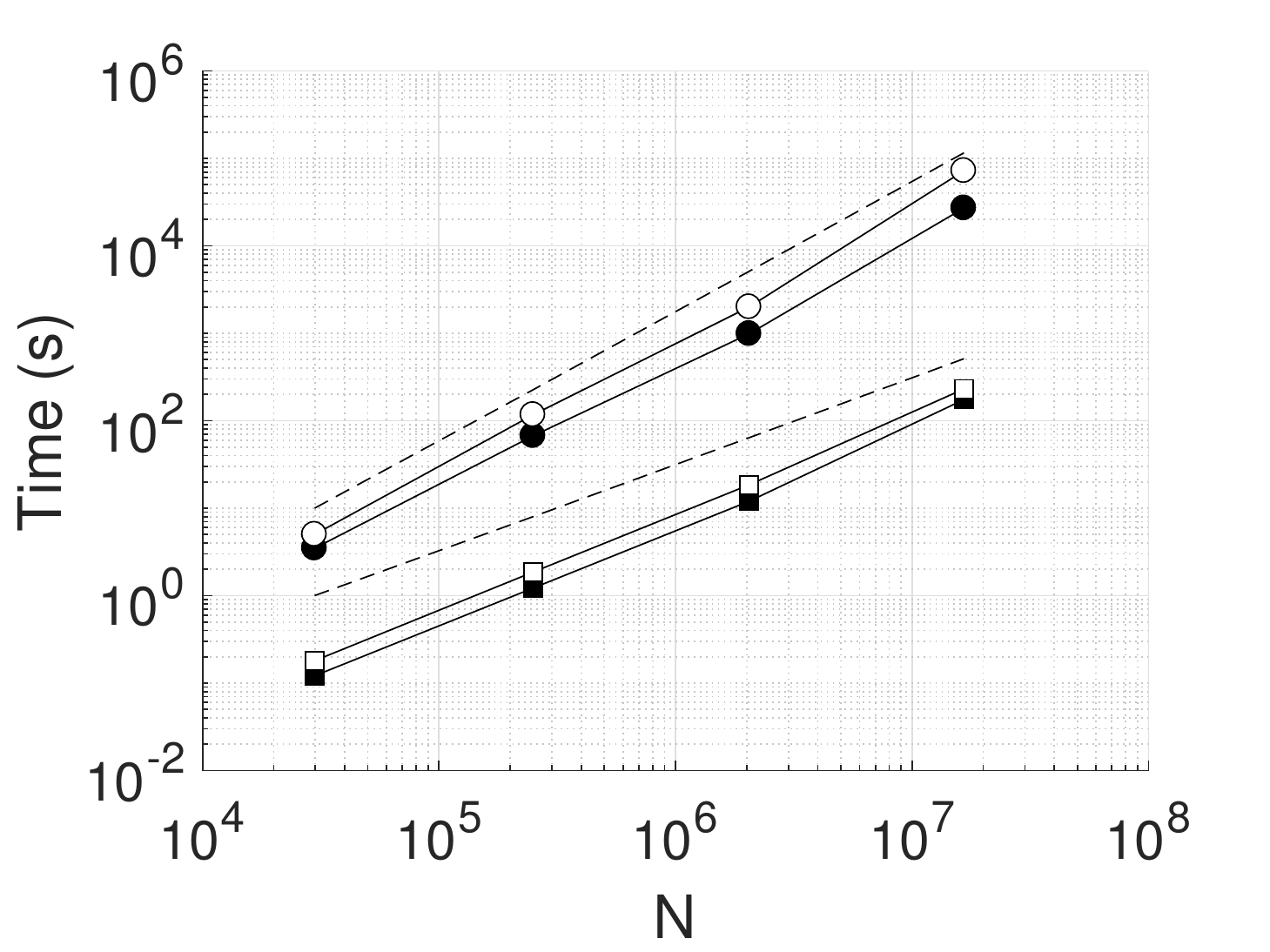}
\end{subfigure}
\hspace{2em}
\begin{subfigure}{0.4\textwidth}
\includegraphics[width=\textwidth]{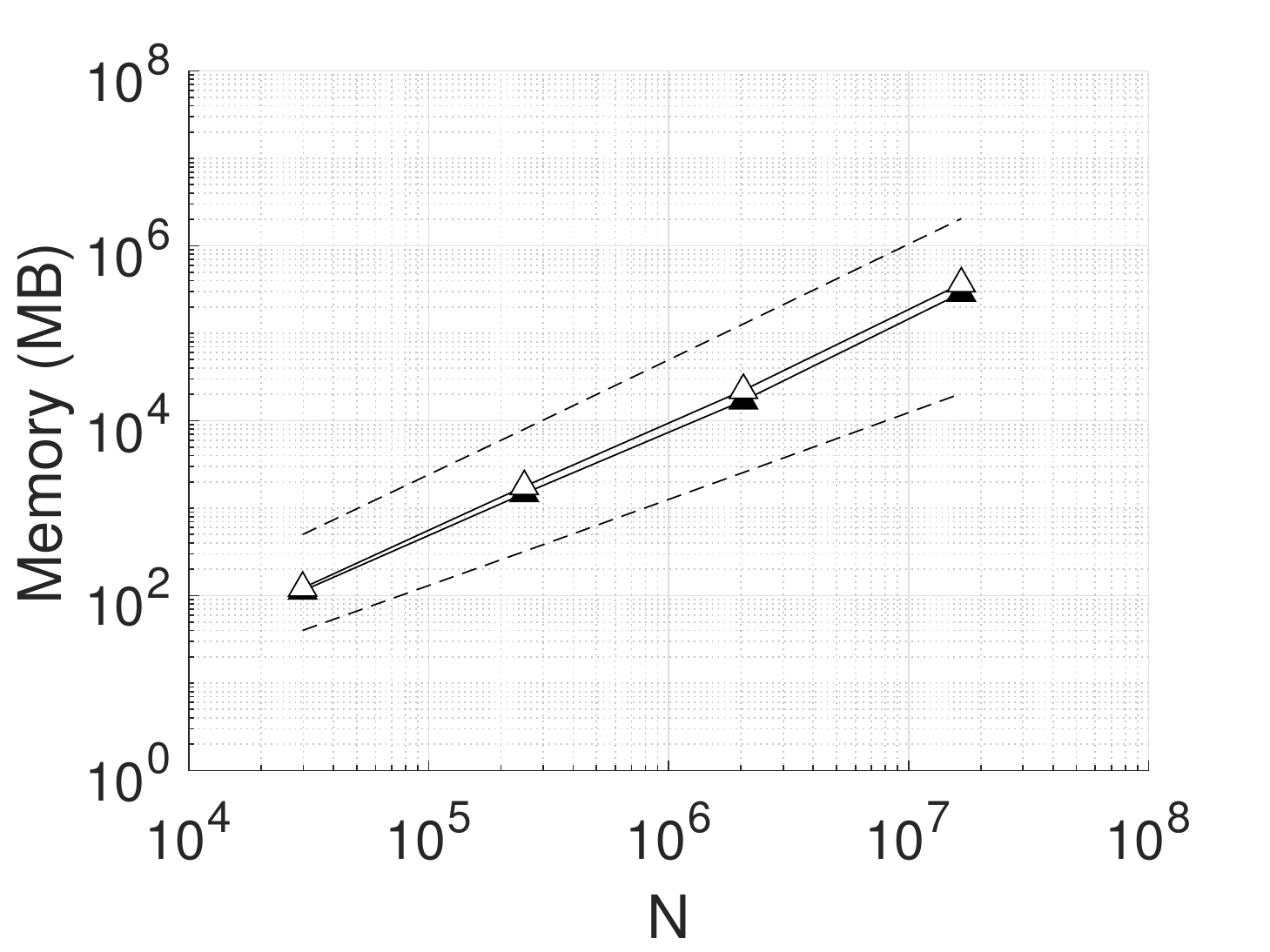}
\end{subfigure}
\caption{Scaling results for 3D example at $\epsilon=10^{-6}$; notation as in Figure \ref{2dcomplexity}. Reference lines (dashed) include $O(N)$ and $O(N^{3/2})$ for timings, and $O(N)$ and $O(N^{4/3})$ for storage.}
\label{3dcomplexity}
\end{figure}

Now consider the analogous problem in 3D. Numerical results are shown in Table \ref{t_results3d} with scaling plots in Figure \ref{3dcomplexity}. The overall conclusions are very similar as before, with PHIF yielding smaller $e_s$ and $n_i$, especially at low accuracy. We highlight in particular the comparison between PHIF at $\epsilon = 10^{-2}$ and HIF at $\epsilon = 10^{-6}$, in which PHIF achieves comparable $n_i$ despite significantly fewer skeletons $|\S_L|$ and worse $e_s$. Furthermore, note the rapid growth of $e_s$ with increasing $N$ for HIF; this is suppressed to a large extent in PHIF, hence producing effective preconditioners that maintain essentially constant $n_i$.

\begin{table}[htb]
\centering
\caption{Numerical results for 3D example.}
\label{t_results3d}
\begin{tabular}{cc|cccc|cccc}
\toprule
&& \multicolumn{4}{|c|}{HIF} & \multicolumn{4}{|c}{PHIF}\\
$\epsilon$ & $N$ & $|\S_L|$ & $e_a$ & $e_s$ & $n_i$ & $|\S_L|$ & $e_a$ & $e_s$ & $n_i$\\
\midrule
\multirow{4}{*}{$10^{-2}$} & $31^3$ & $638$ & $1.7$e$-2$ & $1.0$e$+0$ & $25$ & $773$ & $3.5$e$-3$ & $8.4$e$-1$ & $9$\\
& $63^3$ & $1324$ & $1.9$e$-2$ & $1.0$e$+0$ & $43$ & $1716$ & $5.7$e$-3$ & $1.0$e$+0$ & $14$\\
& $127^3$ & $2605$ & $1.8$e$-2$ & $9.9$e$-1$ & $89$ & $3585$ & $8.1$e$-3$ & $1.0$e$+0$ & $26$\\
& $255^3$ & $4477$ & $1.7$e$-2$ & $1.0$e$+0$ & $203$ & $7026$ & $7.5$e$-3$ & $9.9$e$-1$ & $43$\\
\midrule
\multirow{4}{*}{$10^{-6}$} & $31^3$ & $1422$ & $9.3$e$-7$ & $6.1$e$-3$ & $10$ & $1573$ & $5.5$e$-7$ & $1.2$e$-4$ & $4$\\
& $63^3$ & $3235$ & $2.8$e$-6$ & $5.4$e$-2$ & $13$ & $3775$ & $1.7$e$-6$ & $2.8$e$-4$& $3$\\
& $127^3$ & $6809$ & $1.4$e$-5$ & $4.3$e$-1$ & $16$ & $8792$ & $2.8$e$-6$ & $5.6$e$-4$ & $3$\\
& $255^3$ & $13726$ & $2.8$e$-5$ & $7.1$e$-1$ & $18$ & $18412$ & $2.3$e$-6$ & $1.3$e$-3$ & $3$\\
\midrule
\multirow{4}{*}{$10^{-10}$} & $31^3$ & $1967$ & $1.2$e$-10$ & $6.7$e$-6$ & $2$ & $2169$ & $8.6$e$-11$ & $1.3$e$-8$ & $3$\\
& $63^3$ & $5112$ & $2.5$e$-10$ & $2.2$e$-5$ & $5$ & $5684$ & $1.4$e$-10$ & $1.7$e$-8$ & $3$\\
& $127^3$ & $12559$ & $7.2$e$-10$ & $1.1$e$-4$ & $8$ & $14204$ & $3.0$e$-10$ & $3.0$e$-8$ & $3$\\
& $255^3$ & $25968$ & $1.2$e$-9$ & $1.3$e$-4$ & $8$ & $30946$ & $4.0$e$-10$ & $3.8$e$-8$ & $3$\\
\bottomrule
\end{tabular}
\end{table}

The data certify that $|\S_L|$ scales as $O(N^{1/3})$, corresponding to formal $O(N \log N)$ complexity. However, we find only an empirical scaling of roughly $O(N^{1.4})$ for both HIF and PHIF. This is most likely due to non-asymptotic effects; at any rate, PHIF does not appear to incur any additional asymptotic cost.

\section{Conclusions}\label{conclusions}
In this paper, we have presented PHIF, a recursively preconditioned version of HIF based on simply adding a block Jacobi preconditioning step before each level of skeletonization. This leads to dramatic improvements in the solve accuracy and therefore in its effectiveness as a direct solver or preconditioner, especially for ill-conditioned problems. Importantly, it retains the near-optimal computational complexity of HIF; this makes it the first of what we anticipate will become a highly successful emerging class of fast structured ND methods with enhanced robustness to the underlying system conditioning.

It is worth emphasizing that PHIF is much more than just applying HIF to a standard ``one-level'' preconditioned matrix. Rather, it seeks to control the condition number at all levels, similarly to \cite{Agullo,Owhadi,Xia:2017,Xing}. Although this control is presently not as tight as in \cite{Agullo,Xia:2017,Xing}, recall that we are addressing a fundamentally different problem with higher-dimensional geometric effects. Future work may be able to bridge this gap.

The PHIF framework is quite general and can be modified in various ways. For instance, any local
preconditioner, e.g., diagonal, can be used in place of block Cholesky; this can help to save on
computational costs, as discussed in Section \ref{Numerical}. Furthermore, the algorithm can, in
principle, extend to non-symmetric and indefinite matrices via local LU factorizations, though some
subtleties may inevitably arise. More research is required to fully settle such matters and also to
explore other important extensions including to structured dense matrices. On this latter point in
particular, fast 2D and 3D integral equation solvers like \cite{HIFIE,Minden} are based on the same
style of multilevel skeletonization; to what extent can the same preconditioning ideas apply? The
principal roadblock would appear to be the cost of modifying the far field; a clever trick to bypass
this could be of great significance.

\section*{Acknowledgements}
The work of J.F. is partially supported by "la Caixa" Fellowship, LCF/BQ/AA16/11580045, sponsored by "la Caixa" Banking
Foundation. The work of L.Y. is partially supported by U.S. Department of Energy, Office of
Science, Office of Advanced Scientific Computing Research, Scientific Discovery through Advanced
Computing (SciDAC) program and the National Science Foundation under award DMS-1818449.


\end{document}